# Skeleton-stabilized ImmersoGeometric Analysis for incompressible viscous flow problems

Tuong Hoang[a,b,e,*], Clemens V. Verhoosel[a], Chao-Zhong Qin[a],
Ferdinando Auricchio[c], Alessandro Reali[c,d], E. Harald van Brummelen[a]

[a]*Eindhoven University of Technology – Department of Mechanical Engineering,
P.O. Box 513, 5600MB Eindhoven, The Netherlands*

[b]*IUSS – Istituto Universitario di Studi Superiori Pavia, 27100 Pavia, Italy*

[c]*University of Pavia – Department of Civil Engineering and Architecture, 27100 Pavia, Italy*

[d]*Technische Universität München – Institute for Advanced Study, 85748 Garching, Germany*

[e]*Utrecht University – Department of Earth Sciences,
P.O. Box 80.021, 3508TA Utrecht, The Netherlands*

**Abstract**

A Skeleton-stabilized ImmersoGeometric Analysis technique is proposed for incompressible viscous flow problems with moderate Reynolds number. The proposed formulation fits within the framework of the finite cell method, where essential boundary conditions are imposed weakly using a Nitsche-type method. The key idea of the proposed formulation is to stabilize the jumps of high-order derivatives of variables over the skeleton of the background mesh. The formulation allows the use of identical finite-dimensional spaces for the approximation of the pressure and velocity fields in immersed domains. The stability issues observed for inf-sup stable discretizations of immersed incompressible flow problems are avoided with this formulation. For B-spline basis functions of degree $k$ with highest regularity, only the derivative of order $k$ has to be controlled, which requires specification of only a single stabilization parameter for the pressure field. The Stokes and Navier-Stokes equations are studied numerically in two and three dimensions using various immersed test cases. Oscillation-free solutions and high-order optimal convergence rates can be obtained. The formulation is shown to be stable even in limit cases where almost every elements of the physical domain is cut, and hence it does not require the existence of interior cells. In terms of the sparsity pattern, the algebraic system has a considerably smaller stencil than counterpart approaches based on Lagrange basis functions. This important property makes the proposed skeleton-stabilized technique computationally practical. To demonstrate the stability and robustness of the method, we perform a simulation of fluid flow through a porous medium, of which the geometry is directly extracted from 3D $\mu CT$ scan data.

[*]Corresponding author.
Email addresses:
tuong.hoang@iusspavia.it, t.hoang@tue.nl, t.hoang@uu.nl (T. Hoang)
c.v.verhoosel@tue.nl (C.V. Verhoosel),
c.qin@tue.nl (C. Qin),
auricchio@unipv.it (F. Auricchio),
alessandro.reali@unipv.it (A. Reali)
e.h.v.brummelen@tue.nl (E.H. van Brummelen)

# 1. Introduction

Finite Element Analysis (FEA) of incompressible flow problems has been an active topic of research over the last decades, with research interests ranging from theoretical aspects to engineering applications. In recent years, IsoGeometric Analysis (IGA) – a spline-based finite element simulation paradigm proposed by Hughes *et al.* [1] with the original aim of establishing a better integration between Computer-Aided Design (CAD) and FEA – has been studied in the context of incompressible flow problems. Isogeometric analysis of mixed formulations for incompressible flow problems based on *inf-sup* stable velocity-pressure pairs has been studied in detail in the literature, which has led to the development of a range of isogeometric element families, namely: Taylor-Hood elements [2–4], Sub-grid elements [4, 5], H(div)-conforming elements [3, 6–8], and Nédélec elements [3]. These element families have been demonstrated to be suitable for the discretization of incompressible flow problems, by virtue of the fact that they leverage the advantageous mathematical properties of the spline basis functions used in isogeometric analysis [9].

The Finite Cell Method (FCM) – an immersed finite element method introduced by Rank *et al.* [10] – has been found to be a natural companion to isogeometric analysis. The key idea of the FCM is to embed a geometrically complex physical domain of interest into a geometrically simple embedding domain, on which a regular mesh can be built easily. The framework in which IGA and FCM are integrated – first considered by Schillinger, Rank *et al.* [11–13] – is also referred to as immersogeometric analysis [14, 15]. On the one hand immersogeometric analysis facilitates consideration of CAD trimming curves in the context of isogeometric analysis. On the other hand, it enables the construction of high-regularity spline spaces over geometrically and topologically complex volumetric domains, for which analysis-suitable spline parameterizations are generally not available. The isogeometric finite cell method has been applied to various problems in solid and structural mechanics (see [16, 17] for comprehensive reviews), in image-based analysis [18, 19], in fluid-structure interaction problems [14, 15, 20], and in various other application areas.

In Hoang *et.al* [21] we have found that when the inf-sup stable isogeometric element families for incompressible flow problems are applied in the finite cell setting, local pressure oscillations generally occur in the vicinity of cut boundaries. An illustration of this oscillatory behavior is shown in Figure 1a. When employing the Galerkin-Least square (GLS) method, we observe similar behavior, as shown in Figure 1b. The occurrence of such oscillations on cut elements with relatively large volume fractions implies that this problem is related to the inf-sup stability of the discrete problem, rather than to conditioning issues related to cut elements with small volume fractions. It is important to note that although the inf-sup stable discretization pairs (and GLS) lead to close-to-optimal converge behavior of global error measures, the oscillations in the pressure field near the immersed boundaries persist under mesh refinement. Hence, the approximation of quantities of interest related to the immersed boundaries is below standard. Moreover, for complex topological domains, *e.g.*, porous media, for which immersed methods can be attractive, cut elements can appear almost everywhere. Consequently, cell-based techniques such as inf-sup stable spaces and GLS-type stabilizations can lead to global unphysical spurious pressure solutions, thereby prohibiting successful application to a large class of immersed incompressible flow problems.

In this manuscript we propose an alternative formulation – based on the skeleton-based stabilization technique developed by Hoang *et.al* [22, 23] in the context of conforming isogeometric analysis – to resolve the stability problems associated with immersed inf-sup isogeometric discretization pairs. In this formulation – which can be regarded as a high-regularity generalization of the continuous interior penalty method by Burman and Hansbo [24] – we rely on stabilization of the mixed form problem by amending the formulation with a skeleton-based penalty term. This alternative form of stabilization relaxes the compatibility constraints on the function spaces to be used, which allows us to consider identical discretization spaces for both the velocity and the pressure fields. Our work is related to developments that have been made in the context of XFEM and CutFEM [25–31]. A novelty of our work is that we fully exploit the maximum regularity of the B-spline basis functions used in IGA, as a consequence of which the introduced stabilization operator only acts on the interface jumps of the highest order derivative of the basis functions. As a result, we only require a single stabilization term to control both the inf-sup stability of the mixed problem and its related pressure-space conditioning issues. We herein propose to supplement the skeleton-based stabilized formulation in [22] with a ghost-penalty term for the velocity space, which is not required from the inf-sup condition point of view, but which is essential to control the robustness of the technique for general cut-cell configurations, and also to control the conditioning of the discretized problem.

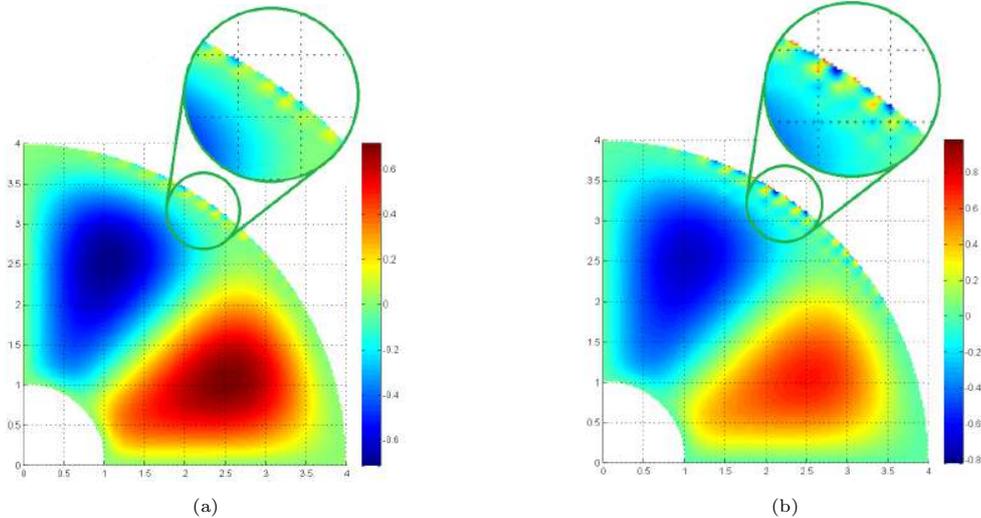

(a)                                         (b)

Figure 1: Unphysical pressure oscillations are observed when the Stokes problem is solved in the standard finite cell setting using *(a)* inf-sup stable isogeometric elements, shown here for Taylor-Hood (see Ref. [21] for details), and *(b)* using a Galerkin Least-squares method.

This paper is outlined as follows. In Section 2 we commence with the introduction of the steady incompressible Navier-Stokes equations and the essentials of the finite cell method. In Section 3 we then present the skeleton-stabilized formulation developed and studied in this work. In Section 4 we discuss the algebraic form of of the developed stabilized formulation, and its effect on the sparsity structure of the system to be solved. The proposed formulation is studied by a series of numerical test cases in Section 5, including the case of a three-dimensional image-based analysis of a microstructural porous medium flow. Conclusions are finally presented in Section 6.

## 2. Preliminaries

Before we introduce the skeleton-based stabilized formulation in Section 3, we herein state the problem setting for the steady Navier-Stokes equations, and the fundamental concepts of the isogeometric finite cell method.

### 2.1. The steady incompressible Navier-Stokes equations

We consider the steady incompressible Navier-Stokes equations on the open bounded domain $\Omega \in \mathbb{R}^d$, where $d = 2, 3$ denotes the spatial dimension of the domain. The Lipschitz boundary $\partial \Omega$ is split into the Dirichlet boundary, $\Gamma_D$, and the Neumann boundary, $\Gamma_N$, such that $\overline{\Gamma_D} \cup \overline{\Gamma_N} = \partial \Omega$ and $\Gamma_D \cap \Gamma_N = \emptyset$. The unit normal vector to $\partial \Omega$, which points out of the domain, is denoted by $\mathbf{n}$. The Navier-Stokes equations for the velocity field $\mathbf{u} : \Omega \to \mathbb{R}^d$ and pressure field $p : \Omega \to \mathbb{R}$ read:

$$\begin{cases} \text{Find } \mathbf{u} : \Omega \to \mathbb{R}^d, \text{ and } p : \Omega \to \mathbb{R} \text{ such that:} \\ \begin{aligned} \nabla \cdot (\mathbf{u} \otimes \mathbf{u}) - \nabla \cdot (2\mu \nabla^s \mathbf{u}) + \nabla p &= \mathbf{f} & \text{in } \Omega, \\ \nabla \cdot \mathbf{u} &= 0 & \text{in } \Omega, \\ \mathbf{u} &= \mathbf{g} & \text{on } \Gamma_D, \\ 2\mu \nabla^s \mathbf{u} \cdot \mathbf{n} - p\mathbf{n} &= \mathbf{h} & \text{on } \Gamma_N. \end{aligned} \end{cases} \quad (1)$$

In this problem formulation the symmetric gradient of the velocity field is denoted by $\nabla^s \mathbf{u} := \frac{1}{2}\left(\nabla \mathbf{u} + (\nabla \mathbf{u})^T\right)$ and $\mu$ represents the kinematic viscosity. The exogenous data $\mathbf{f} : \Omega \to \mathbb{R}^d$, $\mathbf{g} : \Gamma_D \to \mathbb{R}^d$, and $\mathbf{h} : \Gamma_N \to \mathbb{R}^d$ represent body forces, prescribed velocity, and traction data, respectively.

To provide a framework for the derivation of the immersed formulation introduced in the next section, we first present the weak formulation in the conforming setting. The weak formulation of the boundary value problem (1) then follows as:

$$\begin{cases} \text{Find } \mathbf{u} \in \boldsymbol{\mathcal{V}}_{\mathbf{g},\Gamma_D} \text{ and } p \in \mathcal{Q}, \text{ such that:} \\ \begin{aligned} c(\mathbf{u}; \mathbf{u}, \mathbf{w}) + a(\mathbf{u}, \mathbf{w}) + b(p, \mathbf{w}) &= \ell_1(\mathbf{w}) & \forall \mathbf{w} \in \boldsymbol{\mathcal{V}}_{\mathbf{0},\Gamma_D}, \\ b(q, \mathbf{u}) &= 0 & \forall q \in \mathcal{Q}. \end{aligned} \end{cases} \quad (2)$$

The linear operators in this formulation are defined as

$$c(\mathbf{v}; \mathbf{u}, \mathbf{w}) := (\mathbf{v} \cdot \nabla \mathbf{u}, \mathbf{w}), \quad (3a)$$

$$a(\mathbf{u}, \mathbf{w}) := 2\mu \left(\nabla^s \mathbf{u}, \nabla^s \mathbf{w}\right), \quad (3b)$$

$$b(q, \mathbf{w}) := -\left(q, \text{div}\, \mathbf{w}\right), \quad (3c)$$

$$\ell_1(\mathbf{w}) := (\mathbf{f}, \mathbf{w}) + \langle \mathbf{h}, \mathbf{w}\rangle_{\Gamma_N}, \quad (3d)$$

where $(\cdot, \cdot)$ denotes the inner product in $L^2(\Omega)$ and $\langle \cdot, \cdot \rangle_{\Gamma_N}$ denotes the inner product in $L^2(\Gamma_N)$. The function spaces in (2) are defined as

$$\boldsymbol{\mathcal{V}}_{\mathbf{g},\Gamma_D} := \left\{\mathbf{u} \in [H^1(\Omega)]^d : \mathbf{u} = \mathbf{g} \text{ on } \Gamma_D\right\}, \qquad \mathcal{Q} := L^2(\Omega), \quad (4)$$

and the velocity test space, $\boldsymbol{\mathcal{V}}_{\mathbf{0},\Gamma_D}$, is taken as the homogeneous version of $\boldsymbol{\mathcal{V}}_{\mathbf{g},\Gamma_D}$. In the case of pure Dirichlet boundary conditions the pressure is determined up to a constant, which then requires supplementation of the additional pressure condition:

$$\mathcal{Q} := L^2_0(\Omega) \equiv \left\{ q \in L^2(\Omega) : \int_\Omega q \, d\Omega = 0 \right\}. \quad (5)$$

## 2.2. The finite cell method

In the finite cell method, the physical domain of interest, $\Omega$, is immersed into a geometrically simple ambient domain, $\mathcal{A} \supset \Omega$, as illustrated in Figure 2a. In this manuscript we consider the ambient domain to be rectangular in 2D or box-shaped in 3D, so that it can be partitioned by a regular grid with uniform spacing $h > 0$ (but there is no restriction on mapping the ambient domain to more complex geometries). We refer to this partitioning as the ambient domain mesh, $\mathcal{T}^h_\mathcal{A}$. Elements in this ambient domain mesh that do not intersect the physical domain are discarded in the finite cell analysis, which leads to the definition of the finite cell background mesh:

$$\mathcal{T}^h := \{K \in \mathcal{T}^h_\mathcal{A} : K \cap \Omega \neq \emptyset\} \quad (6)$$

The ambient domain mesh and background mesh are illustrated in Figure 2.

The conceptual idea of the finite cell method is to construct a suitable discretization space on the background mesh, and to use that basis in a Galerkin formulation pertaining to the physical domain. Dirichlet boundary conditions on non-conforming edges are typically enforced weakly, most commonly by means of Nitsche's method [32]. We will introduce the Nitsche formulation for the problem (1) in Section 3. In the remainder of this section we introduce the B-spline basis defined over the background mesh, and the integration procedure employed to evaluate volume and (immersed) surface integrals over elements that are cut by the immersed boundary.

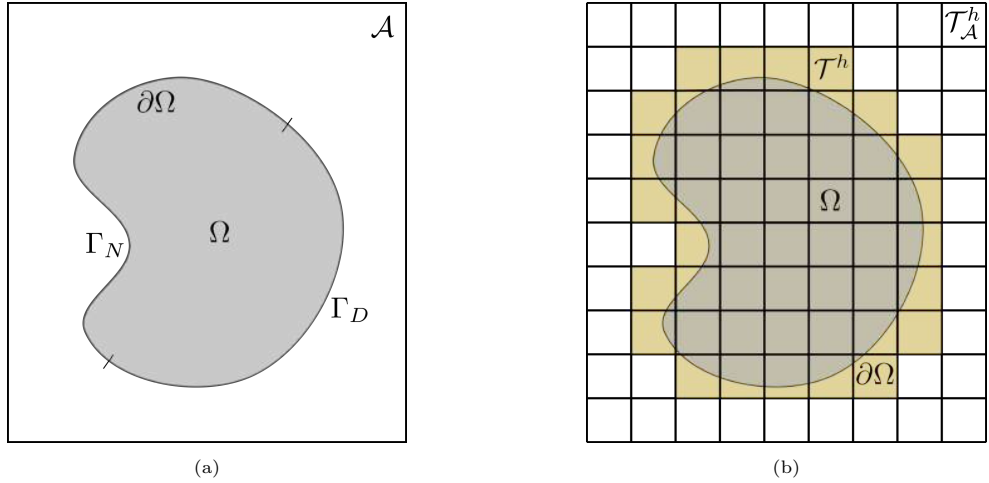

Figure 2: Schematic representation of *(a)* the physical domain $\Omega$ and ambient domain $\mathcal{A}$ as considered in the finite cell method, and *(b)* the ambient domain mesh, $\mathcal{T}_\mathcal{A}^h$ (covering the complete ambient domain), and background mesh, $\mathcal{T}^h$ (marked by the yellow background shading).

*2.2.1. B-spline basis*

By virtue of the fact that in the finite cell method basis functions are constructed on a regular background mesh, it enables the isogeometric analysis of complex-shaped physical domains. In this manuscript we restrict ourselves to a single patch open B-spline basis over the ambient domain mesh, defined by non-decreasing knot vectors in all spatial directions $\delta = 1, \ldots, d$,

$$\Xi^\delta = [\,\underbrace{\xi_1^\delta, \ldots, \xi_1^\delta}_{(k+1)-\text{times}}, \xi_2^\delta, \ldots, \xi_{m^\delta-1}^\delta, \underbrace{\xi_{m^\delta}^\delta, \ldots, \xi_{m^\delta}^\delta}_{(k+1)-\text{times}}], \tag{7}$$

In accordance with the definition of open B-splines the first and last knot values are repeated $k+1$ times, where $k$ denotes the global isotropic polynomial degree of the basis. As explained above, we herein decide to align the knot vectors with the ambient domain, which essentially implies that the we have an identity geometric map between the parameter domain and the ambient domain. The spacing between two consecutive knots is therefore equal to the global isotropic mesh parameter $h$.

Using the knot vectors (7) a B-spline basis of degree $k$ can be constructed over the ambient domain by means of the recursive Cox-de Boor formula [33]. We denote this B-spline basis by $\mathcal{N}_\mathcal{A}^k = \{N_{\mathcal{A},I}^k : \mathcal{A} \to \mathbb{R}\}_{I=1}^{n_\mathcal{A}}$, where the total number of basis functions is equal to $n_\mathcal{A} = \otimes_{\delta=1}^d \{m^\delta + k - 1\}$, with $m^\delta$ the number of unique knot values per direction. In the finite cell analysis we discard the basis function that are not supported on the background mesh $\mathcal{T}^h$, so that the B-spline basis follows as:

$$\mathcal{N}^k := \{N \in \mathcal{N}_\mathcal{A}^k : \text{supp}(N) \cap \mathcal{T}^h \neq \emptyset\} \tag{8}$$

Note that, by definition, all basis functions in $\mathcal{N}^k$ have positive support over the physical domain $\Omega$. The cardinality of $\mathcal{N}^k$ is denoted by $n \leq n_\mathcal{A}$. We herein consider maximum regularity B-spline bases – as indicated by the non-repeated internal knot values in (7) – so that the basis functions are $C^{k-1}$ continuous. The B-spline function space $\mathcal{S}_h^k$ spanned by the basis $\mathcal{N}^k$ is therefore a finite dimensional subspace of the Sobolev space $H^k(\Omega^h)$, where $\Omega^h = int(\bigcup_{K \in \mathcal{T}^h} \overline{K})$.

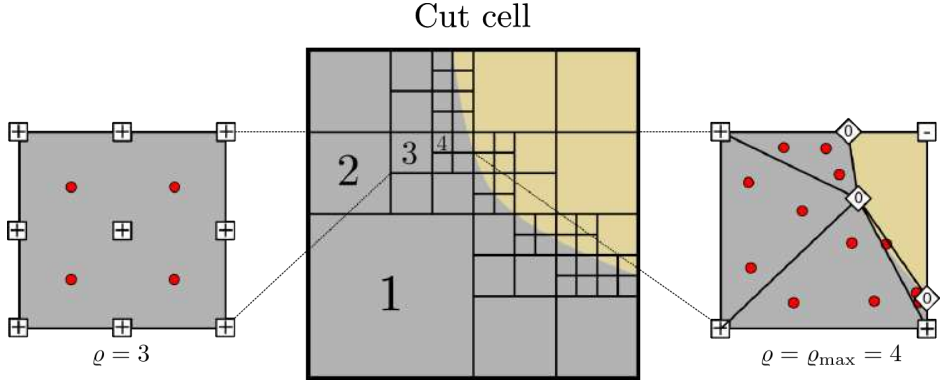

Figure 3: Illustration of the bisection-based tessellation scheme used to generate quadrature rules for the cells that are cut by the immersed boundary. The vertex markers indicate whether the employed interpolant of the level set function is zero (on the boundary), positive (inside the domain), or negative (outside the domain). The red points in the sub-cell zooms are an illustration of the distribution of the integration points in such cells.

*2.2.2. Cut cell integration*

For elements in the background mesh that are intersected by the boundaries of the physical domain standard quadrature rules are inaccurate, since effectively discontinuous functions are integrated over such cut cells. The FCM therefore generally employs an advanced numerical integration technique for cut cells. Herein we use the bisectioning-based segmentation scheme proposed in [19] in the context of the isogeometric finite cell analysis of image-based geometric models, which also enables us to extract a parametrization and quadrature rules for the immersed boundaries.

We illustrate the bisection-based tessellation scheme in Figure 3 for completeness. The element-by-element routine commences with the evaluation of a level set function in the vertices associated with a $\varrho_{\max}$-times uniform refinement of the element. This level set can either be derived from voxel data in a scan-based analysis, or a signed-distance function can be considered in the case that the geometry is provided by a CAD model. The integration points for a cut cell are assembled by traversing the levels of uniform refinement, where for each sub-cell in that level it is determined whether the interface passes through it. If all vertices of a sub-cell exceed a specified threshold value (zero in the case of a signed-distance function) the sub-cell is kept as an integration sub-cell. Otherwise a further subdivision of the sub-cell is considered, and the same check is performed on the next level. On the lowest level this recursion is closed with a tessellation procedure. From an implementation perspective the integration points and weights on all sub-cells are collected on the level of the cut element, which essentially provides us with an integration scheme tailored to the cut element.

The tessellation procedure used on the deepest level of integration refinement enables the extraction of a parametrization of the boundary. In essence, the immersed boundary is reconstructed on an element-by-element basis by identifying the faces of the integration sub-cells that coincide with the immersed boundary. The collection of sub-cell faces that approximates the immersed boundary provides a piece-wise linear parametrization of this boundary. Using this piece-wise parametrization, quadrature rules can be constructed. Evidently, the refinement parameter $\varrho_{\max}$ controls the accuracy with which the geometry is approximated, as well as the computational complexity of the approximation.

## 3. Skeleton-stabilized immersed isogeometric analysis for the Navier-Stokes equations

In this section we introduce the skeleton-stabilized immersed isogeometric analysis formulation for the Navier-Stokes equations. We commence with the definition of the topological structures on which this stabilization technique is based, after which we present the two stabilization aspects in our formulation, *viz.* the ghost penalty stabilization of the velocity components at the cut boundaries, and the pressure stabilization on the skeleton of the background mesh.

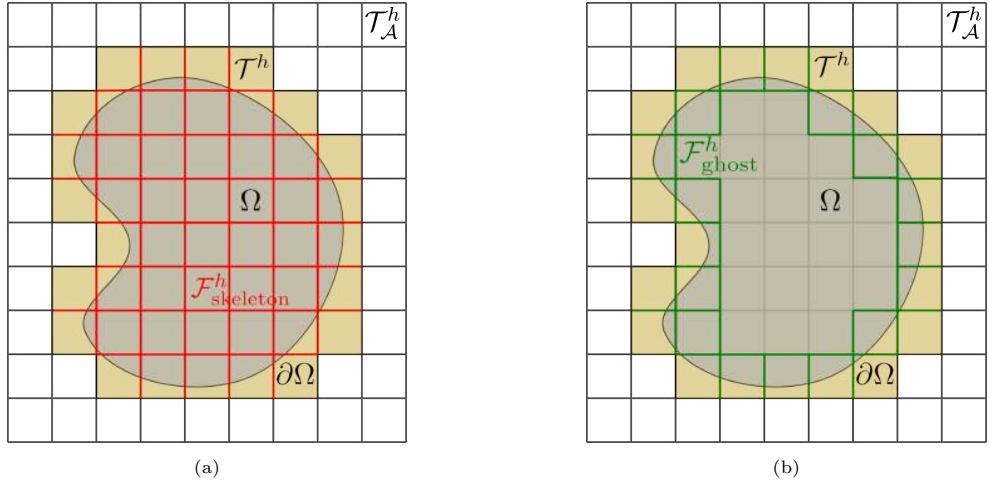

Figure 4: Schematic representation of *(a)* the skeleton structure (red), and *(b)* the ghost interface structure (green).

*3.1. The background mesh: skeleton structure and ghost structure*

We consider the background mesh $\mathcal{T}^h$ as defined in Section 2.2. The stabilized formulation presented herein is based on the skeleton of this background mesh, which is defined as

$$\mathcal{F}^h_{\text{skeleton}} = \{\partial K \cap \partial K' \mid K, K' \in \mathcal{T}^h, K \neq K'\}. \tag{9}$$

This skeleton, for which the mesh parameter $h$ is associated with that of the background mesh, is illustrated in Figure 4a. Note that the boundary faces of the mesh $\mathcal{T}^h$ are not a part of this skeleton mesh.

Besides the skeleton structure (9) we also consider that part of the skeleton which coincides with the faces of all cut cells in the domain, to which we refer as the ghost skeleton:

$$\mathcal{F}^h_{\text{ghost}} = \mathcal{F}^h_{\text{skeleton}} \cap \{F \mid F \subset \partial K : K \in \mathcal{T}^h, K \cap \partial\Omega \neq \varnothing\} \tag{10}$$

This ghost skeleton structure is illustrated in Figure 4b, where the mesh parameter $h$ is still associated with that of the background mesh. Note that this structure can also contain some faces that do not intersect the boundary of the physical domain.

Since we herein consider single patch discretizations of the ambient domain with maximal regularity, the basis functions are $C^{k-1}$ continuous over all faces in the skeleton mesh $\mathcal{F}^h_{\text{skeleton}}$. Therefore, the jumps in the normal derivative up to order $k-1$ of all functions in $\mathcal{S}^k_h = \text{span}(\mathcal{N}^k)$ vanish:

$$[\![\partial^i_n f]\!]_F = 0, \quad 0 \leq i \leq k-1, \qquad \forall f \in \mathcal{S}^k_h, \tag{11}$$

where the jump operator $[\![\cdot]\!]$ associates to any function $f$ in the broken Sobolev space $\{f \in L^2(\Omega) : f|_K \in H^1(K), \forall K \in \mathcal{T}^h\}$ the $L^2(\mathcal{F}^h_{\text{skeleton}})$-valued function:

$$[\![f]\!] = f^+ - f^-$$

The superscripts $(\cdot)^\pm$ refer to the traces of $f$ on the two opposite sides of each face $F \in \mathcal{F}^h_{\text{skeleton}}$, with an arbitrary allocation of $+$ and $-$.

*3.2. The skeleton-stabilized finite cell formulation*

In this contribution we study the discretization of problem (1) using identical highest smoothness ($C^{k-1}$) spline discretizations of degree $k$ for both the velocity and pressure fields:

$$\boldsymbol{\mathcal{V}}^h := \left[\mathcal{S}^k_h\right]^d, \qquad\qquad \mathcal{Q}^h := \mathcal{S}^k_h. \tag{12}$$

The skeleton-stabilized finite cell formulation for the system (1) reads:

$$\begin{cases} \text{Find } \mathbf{u}^h \in \boldsymbol{\mathcal{V}}^h \text{ and } p^h \in \mathcal{Q}^h, \text{ such that:} \\ c(\mathbf{u}^h; \mathbf{u}^h, \mathbf{w}^h) + a^h(\mathbf{u}^h, \mathbf{w}^h) + s^h_{\text{ghost}}(\mathbf{u}^h, \mathbf{w}^h) + b^h(p^h, \mathbf{w}^h) = \ell_1^h(\mathbf{w}^h) \quad \forall \mathbf{w}^h \in \boldsymbol{\mathcal{V}}^h, \\ \qquad\qquad\qquad\qquad\qquad b^h(q, \mathbf{u}^h) - s^h_{\text{skeleton}}(p^h, q^h) = \ell_2^h(q^h) \quad \forall q^h \in \mathcal{Q}^h, \end{cases} \quad (13)$$

where the linear operators as introduced in the conforming setting in equation (3) are supplemented with additional terms for Nitsche's imposition of the boundary conditions [32]:

$$a^h(\mathbf{u}^h, \mathbf{w}^h) := a(\mathbf{u}^h, \mathbf{w}^h) - 2\mu \left[ \langle \nabla^s \mathbf{u}^h \cdot \mathbf{n}, \mathbf{w}^h \rangle_{\Gamma_D} + \langle \nabla^s \mathbf{w}^h \cdot \mathbf{n}, \mathbf{u}^h \rangle_{\Gamma_D} \right] + \mu \langle \beta h^{-1} \mathbf{u}^h, \mathbf{w}^h \rangle_{\Gamma_D}, \quad (14a)$$

$$b^h(q^h, \mathbf{w}^h) := b(q^h, \mathbf{w}^h) + \langle q^h, \mathbf{w}^h \cdot \mathbf{n} \rangle_{\Gamma_D}, \quad (14b)$$

$$\ell_1^h(\mathbf{w}^h) := \ell_1(\mathbf{w}^h) - 2\mu \langle \nabla^s \mathbf{w}^h \cdot \mathbf{n}, \mathbf{g} \rangle_{\Gamma_D} + \mu \langle \beta h^{-1} \mathbf{g}, \mathbf{w}^h \rangle_{\Gamma_D}, \quad (14c)$$

$$\ell_2^h(q^h) := \langle q^h \mathbf{n} \cdot \mathbf{g} \rangle_{\Gamma_D}, \quad (14d)$$

where $\langle \cdot, \cdot \rangle_{\Gamma_D}$ denotes the inner product in $L^2(\Gamma_D)$. The two stabilization operators that are appended to the weak form (13) are defined as:

$$s^h_{\text{ghost}}(\mathbf{u}^h, \mathbf{w}^h) := \sum_{F \in \mathcal{F}_{\text{ghost}}} \int_F \widetilde{\gamma} \mu h^{2k-1} [\![\partial_n^k \mathbf{u}_h]\!] \cdot [\![\partial_n^k \mathbf{v}_h]\!] \, \mathrm{d}s, \quad (15a)$$

$$s^h_{\text{skeleton}}(p^h, q^h) := \sum_{F \in \mathcal{F}^h_{\text{skeleton}}} \int_F \gamma \mu^{-1} h^{2k+1} [\![\partial_n^k p^h]\!] [\![\partial_n^k q^h]\!] \, \mathrm{d}s, \quad (15b)$$

where $\widetilde{\gamma}$ and $\gamma$ denote certain suitable positive stabilization parameters; see Remark 2 below.

The operator $s^h_{\text{ghost}}(\mathbf{u}^h, \mathbf{w}^h)$ is referred to as the ghost-penalty operator [25, 26]. This term – which penalizes the (non-vanishing) jump in the $k$-th order normal derivative on the ghost skeleton, enables scaling of the Nitsche penalty term by the reciprocal mesh size parameter $h^{-1}$ of the *background mesh*, independent of the cut-element configurations. Without the ghost-penalty term, the Nitsche term would have to be based on the reciprocal of the *cut-element* size to ensure stability. However, this would result in configuration-sensitive stability and severe ill-conditioning in critical cases such as sliver cut configurations [26, 34]. The ghost-penalty stabilization effectively controls the conditioning and maintains the accuracy of the velocity field. The condition numbers then scale as in the case of conforming discretizations, and are independent of the configuration of the cut cells [25, 26, 29].

The operator $s^h_{\text{skeleton}}(p^h, q^h)$ in (13) is referred to as the skeleton-penalty operator, which was developed for conforming isogeometric discretizations of the unsteady incompressible Navier-Stokes equations in [22, 23], and is applied here without any modification to the immersed setting. This term allows to use identical pairs of spaces for the velocity and pressure fields, as defined in equation (12). It should be emphasized that the skeleton structure is defined not only inside the physical domain but in the whole background mesh. In this way, the pressure-stabilization not only ensures inf-sup stability over the complete domain but also resolves the conditioning issue related to the pressure field in the case of pathological cut configurations.

**Remark 1.** *We note that the stabilization parameters $\gamma$ and $\widetilde{\gamma}$, and the viscosity parameter $\mu$, in the operators (15a) and (15b) are kept inside the integrand for the sake of generality. For the simulations considered herein – where we focus on moderate Reynolds numbers flows – these stabilization parameters are defined globally. This global scaling may not extend to, e.g., the more general case of convection-dominated flows.*

**Remark 2.** *The positive parameters $\gamma$ and $\widetilde{\gamma}$ are herein selected in an ad hoc manner, where an important guideline is that they should decrease with increasing regularity. In the full regularity setting considered herein, this implies that they decrease with increasing order of the discretization. See Section 5.1 for a more detailed study.*

**Remark 3.** *The pressure skeleton-stabilization term and the ghost penalty stabilization term scale differently with the mesh size, viz. with $h^{2k+1}$ and $h^{2k-1}$, respectively. This difference stems from the* fact that the velocity field resides in $H^1$, while the pressure resides in $L^2$. Note that although we have restricted ourselves herein to regular meshes with global and isotropic size h, there is no fundamental restriction to the application of the formulation (13) in the context of non-uniformly spaced grids.

**Remark 4.** *An important stability property of our skeleton-stabilized formulation is that, different from other approaches in the literature [5, 25–31, 35–37], we herein do not require assumptions on the existence of interior cells (cells not cut by the boundary). In limit cases where almost every element is cut (i.e., $\mathcal{F}_{\text{ghost}}^h \simeq \mathcal{F}_{\text{skeleton}}^h$), numerical experiments suggest that our formulation remains stable, see Section 5.3.*

The rationale behind the skeleton-stabilized formulation – which for sufficiently smooth velocity and pressure solutions ($\mathbf{u} \in [H^k(\Omega)]^d$ and $p \in H^k(\Omega)$) is consistent with (1) – is that it effectively targets the shortcomings observed using inf-sup stable spaces. The skeleton stabilization operator (15b) is tied to the background domain, in the sense that it is completely independent of the shape and volume fraction of the cut cells it pertains to. As a consequence, the stabilizing effect of the operator does not decrease with decreasing volume fractions. This contrasts the situation in which inf-sup stable pairs are considered, since in that setting the cut cell characteristics have been observed to impact the inf-sup stability [21]. Moreover, the stabilization operator (15b) can be conceived of as a weakly imposed constraint on the highest-order non-vanishing derivative of the pressure field, which essentially means that it controls the smoothness of the extension of the interior domain into the exterior domain. In [22] we have demonstrated that the operator (15b) is related to the least squares minimization problem for the highest-order derivative jumps. Thereby the operator effectively suppresses oscillations in the pressure field near the immersed boundaries.

## 4. The algebraic form

To elucidate the structure of the skeleton-stabilized finite cell formulation and to establish the impact of the stabilization terms on the bandwidth of the resulting system of equations, we here consider the algebraic form of the formulation. Let $\{\mathbf{N}_i\}_{i=1}^{n_u}$ and $\{N_i\}_{i=1}^{n_p}$ denote two sets of B-spline basis functions (as defined in Section 2.2) for the velocity and pressure fields, respectively. The vector-valued velocity basis functions are defined as

$$\mathbf{N}_{i=j+\delta n} = N_j \mathbf{e}_\delta, \qquad\qquad j = 1, \ldots, n \text{ and } \delta = 1, \ldots, d \qquad (16)$$

where $n$ is the number of control points, $d$ the number of spatial dimensions (evidently, $n_u = dn$ and $n_p = n$), and $\mathbf{e}_\delta$ is the unit vector in the direction $\delta$. The basis functions span the discrete velocity and pressure spaces

$$\boldsymbol{\mathcal{V}}^h = \text{span}\{\mathbf{N}_i\}_{i=1}^{n_u}, \qquad\qquad \mathcal{Q}^h = \text{span}\{N_i\}_{i=1}^{n_p}. \qquad (17)$$

The discrete velocity and pressure fields can be expressed as

$$\mathbf{u}^h(\mathbf{x}) = \sum_{i=1}^{n_u} \mathbf{N}_i(\mathbf{x})\hat{u}_i, \qquad\qquad p^h(\mathbf{x}) = \sum_{i=1}^{n_p} N_i(\mathbf{x})\hat{p}_i. \qquad (18)$$

where $\hat{\mathbf{u}} = (\hat{u}_1, \hat{u}_2, \ldots, \hat{u}_{n_u})^T$ and $\hat{\mathbf{p}} = (\hat{p}_1, \hat{p}_2, \ldots, \hat{p}_{n_p})^T$ are vectors of degrees of freedom. In the absence of constraints, the number of velocity-degrees of freedom, $n_u$, is $d$ times that of pressure degrees of freedom, $n_p$.

The formulation (13) can be cast into an algebraic system of equations of size $n = n_u + n_p$:

$$\begin{cases} \text{Find } \hat{\mathbf{u}} \in \mathbb{R}^{n_u} \text{ and } \hat{\mathbf{p}} \in \mathbb{R}^{n_p}, \text{ such that:} \\ [\mathbf{C}(\hat{\mathbf{u}}) + \mathbf{A} + \mathbf{S}_{\text{ghost}}]\hat{\mathbf{u}} + \mathbf{B}^T\hat{\mathbf{p}} = \mathbf{f}_1, \\ \qquad\qquad \mathbf{B}\hat{\mathbf{u}} - \mathbf{S}_{\text{skeleton}}\hat{\mathbf{p}} = \mathbf{f}_2. \end{cases} \qquad (19)$$

We employ Picard iterations to solve this nonlinear algebraic problem. The matrices and vectors in (19) pertaining to the standard volume and boundary surface terms can be expressed in terms of the operators (14) as:

$$C(\hat{\mathbf{u}})_{ij} := c(\hat{\mathbf{u}}; \mathbf{N}_j, \mathbf{N}_i), \tag{20a}$$
$$A_{ij} := a^h(\mathbf{N}_j, \mathbf{N}_i), \tag{20b}$$
$$B_{ij} := b(N_i, \mathbf{N}_j), \tag{20c}$$
$$f_{1,i} := \ell_1^h(\mathbf{N}_i) \tag{20d}$$
$$f_{2,i} := \ell_2(N_i). \tag{20e}$$

The stabilization matrices in (19) pertain to the skeleton and ghost structure of the background mesh, $\mathcal{F}^h_{\text{skeleton}}$ and $\mathcal{F}^h_{\text{ghost}}$, respectively, and hence require data structures to evaluate the jump of high-order derivatives of the basis functions across the background mesh element interfaces through the operators in (15) as:

$$S_{\text{skeleton},ij} = s^h_{\text{skeleton}}(N_j, N_i), \tag{21a}$$
$$S_{\text{ghost},ij} = s^h_{\text{ghost}}(\mathbf{N}_j, \mathbf{N}_i). \tag{21b}$$

Due to the fact that the jump operators on the highest-order derivatives of the B-spline basis functions provide additional connectivity between basis functions, the stabilization matrices (21) have an effect on the sparsity pattern of the algebraic problem. In Figure 5 we present a comparison of the sparsity patterns of the system matrices in two dimensions for the cases of a second-order ($k = 2$) B-spline basis as considered herein and a second-order ($k = 2$) Lagrange basis (closely resembling the continuous interior penalty method). Note that since both bases are constructed over the same background mesh, the number of Lagrange basis functions is significantly larger than the number of B-spline basis functions, by virtue of the fact that, as opposed to Lagrange basis functions, for full-regularity B-splines the number of basis functions does not scale proportionally with the degree of the basis to the power $d$.

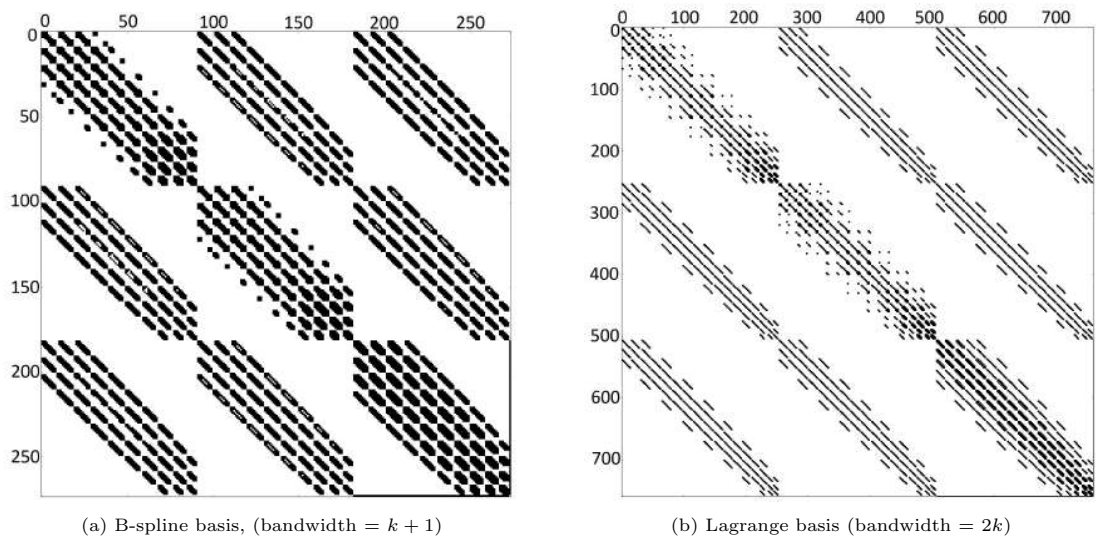

(a) B-spline basis, (bandwidth = $k+1$)     (b) Lagrange basis (bandwidth = $2k$)

Figure 5: Illustration of the sparsity patterns of the system matrix corresponding to (19) for *(a)* a second-order B-spline basis, and *(b)* a second-order Lagrange basis. The stencil of the B-spline case is smaller than that of the Lagrange case. See Figure 6 for details. Note that here both bases are constructed over the same background mesh and therefore result in different numbers of degrees of freedom; the figure sizes are thus not of the same scale.

Inspection of the velocity-velocity and pressure-pressure blocks reveals that the footprint of the stabilization operators have a different effect for the two bases, in the sense that for the $k = 2$ case $2k = 4$ off-diagonal bands are observed for the Lagrange basis, and $k + 1 = 3$ for the B-spline basis. This difference – which becomes more pronounced when the degree $k$ increases – was also observed in the mesh conforming case in Ref. [22], with the difference that in the immersed setting both the velocity and pressure space are stabilized, thereby making the impact of the stabilization operators on the computational effort larger in the immersed setting. Following the arguments in Ref. [22], the difference in number of off-diagonal bands can be explained by comparison of the one-dimensional B-spline and Lagrange bases, as shown for the cubic ($k = 3$) case in Figure 6. This figure corroborates that in the case of full-regularity B-splines all derivative jumps up to order $k - 1$ vanish across element interfaces, as a result of which only the $k$-th derivative jump operator impacts the sparsity pattern. The number of off-diagonal bands for the stabilization operators is therefore in this case equal to $k + 1$. In contrast, for Lagrange bases, the lower-order derivative jumps are non-vanishing, as a result of which $2k$ off-diagonal bands appear.

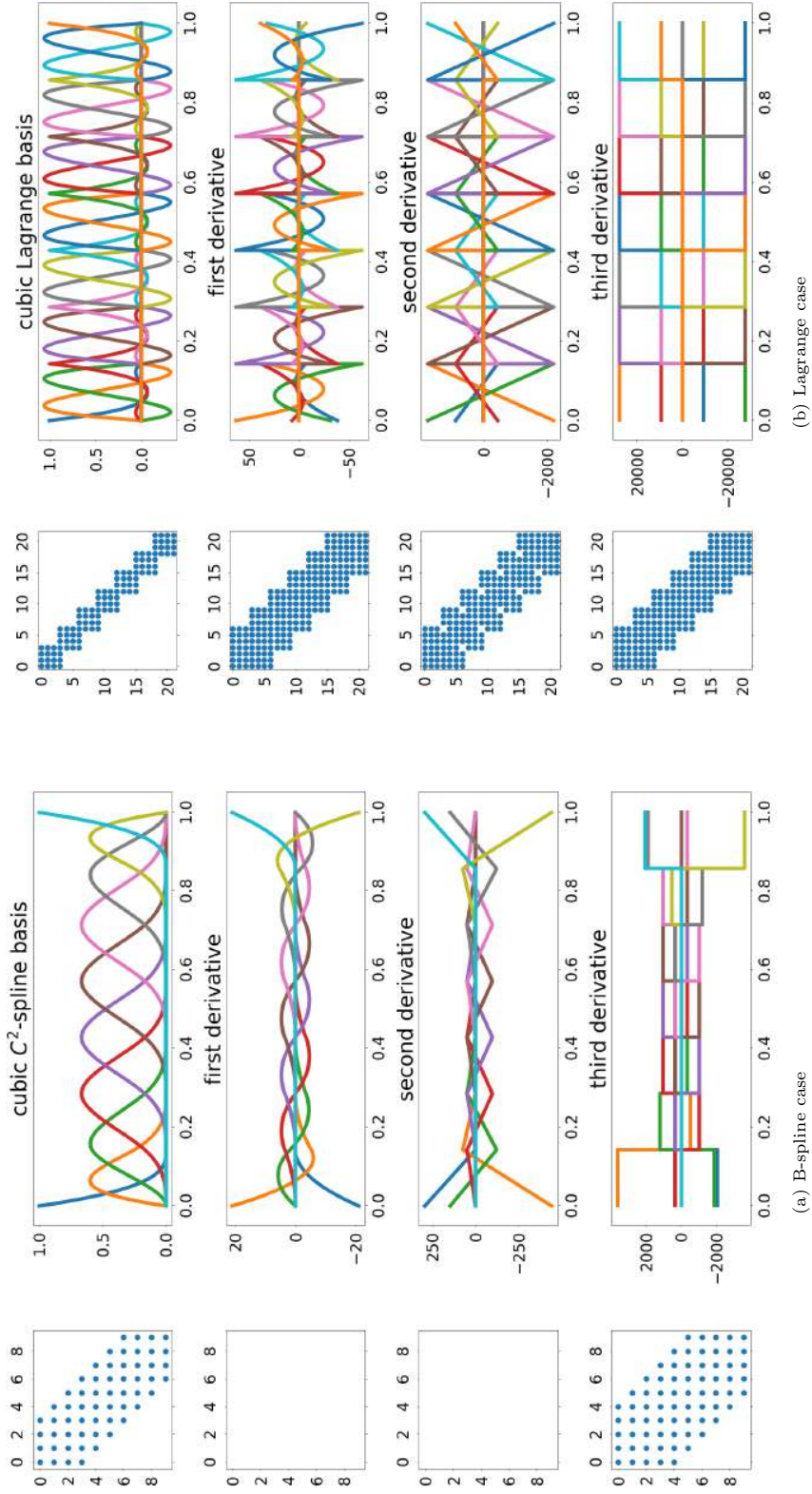

Figure 6: Comparison of the (a) univariate cubic B-spline basis of full regularity, and (b) univariate cubic Lagrange basis. Due to the $C^{k-1}$ continuity of the B-spline basis only the $k$-th derivative jump is non-vanishing, this in contrast to the case of the $C^0$ continuous Lagrange basis. Moreover, the bandwidth in the B-spline case is $k+1$, much smaller than in the Lagrange case which is $2k$. The visualized sparsity patterns are of the stiffness matrix (first row) and of the first, second, and third order jump-derivative matrices (second, third, and last rows, respectively).

## 5. Numerical simulations

In this section we investigate the numerical performance of the proposed immersed skeleton-stabilized formulation. In all cases, the system (13) is solved using identical highest-regularity B-spline spaces for the approximations of the velocity and pressure fields. Unless stated otherwise, the number of bi-sectioning levels that determines the accuracy of the geometry representation (see Section 2.2.2) is taken equal to six. Evidently, when studying higher-order approximations, one ideally wants to resolve the geometry as closely as possible. The above-mentioned bi-sectioning depth is chosen such that the simulations remain computationally tractable.

### 5.1. Steady Navier-Stokes flow in a quarter annulus domain

We consider the steady Navier-Stokes equations on an open quarter-annulus domain

$$\Omega = \left\{ (x, y) \in \mathbb{R}^2_{>0} : R_1^2 < x^2 + y^2 < R_2^2 \right\},$$

with inner radius $R_1 = 1$ and outer radius $R_2 = 4$; see Figure 7. Dirichlet boundary conditions are prescribed on the entire boundary $\partial \Omega = \Gamma_D$. The body force $\mathbf{f}$ and Dirichlet data $\mathbf{g}$ are selected in accordance with the manufactured solution

$$\begin{aligned}
u_1 &= 10^{-6} x^2 y^4 (x^2 + y^2 - 1)(x^2 + y^2 - 16)(5x^4 + 18x^2 y^2 - 85x^2 + 13y^4 - 153y^2 + 80), \\
u_2 &= 10^{-6} xy^5 (x^2 + y^2 - 1)(x^2 + y^2 - 16)(102x^2 + 34y^2 - 10x^4 - 12x^2 y^2 - 2y^4 - 32), \\
p &= 10^{-7} xy(y^2 - x^2)(x^2 + y^2 - 16)^2 (x^2 + y^2 - 1)^2 \exp\left(14(x^2 + y^2)^{-1/2}\right),
\end{aligned} \quad (22)$$

of problem (1) without the inertia term and with viscosity $\mu = 1$. This manufactured solution is adopted from Refs. [21, 38]. Note that $u_1$ and $u_2$ vanish on $\partial \Omega$, and hence $\mathbf{g} = \mathbf{0}$. Moreover, the manufactured pressure solution complies with the zero-average pressure condition $\int_\Omega p = 0$, which is imposed here by means of a Lagrange multiplier.

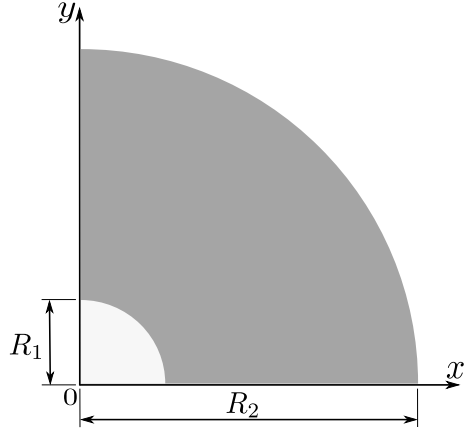

Figure 7: Geometry definition of the quarter-annulus ring problem.

We tested this problem for the Stokes case in Ref. [21] with different families of inf-sup stable isogeometric spaces. A representative result for that setting was shown in Figure 1a, from which unphysical pressure oscillations in the vicinity of the cut elements are clearly observed. Using inf-sup stable pairs, similar oscillations are also observed in the Navier-Stokes case (results not displayed). In contrast, the pressure field computed using the skeleton-stabilized formulation (13) – illustrated in Figure 8 for the case of quadratic B-splines with a $21 \times 21$ elements ambient domain mesh – is free of oscillations. Note that the physical domain is completely immersed in the ambient domain, in the sense that none of the boundaries conform to the background mesh.

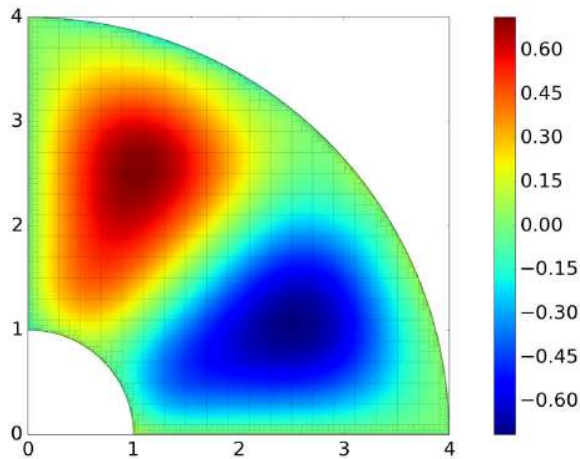

Figure 8: Pressure solution of the steady Navier-Stokes equations with $\mu = 1$, computed using the skeleton-stabilized formulation with quadratic B-splines. The original ambient domain mesh consists of $21 \times 21$ elements.

In Figure 9 we present mesh convergence results for the proposed stabilized formulation, where a sequence of uniformly refined background meshes is generated starting from the $11 \times 11$ elements coarsest ambient domain mesh. The finest level ambient domain mesh contains $176 \times 176$ elements. We consider $k = 1, 2, 3$ full-regularity B-splines with stabilization parameters $\gamma = 10$ for $k = 1$, $\gamma = 0.1$ for $k = 2$, $\gamma = 5 \times 10^{-4}$ for $k = 3$ and $\widetilde{\gamma} = 10^{-k-1}$ for all $k$. These stabilization parameters have been selected based on a numerical sensitivity study (see details below). We observe optimal rates of convergence of $k + 1$ and $k$ for the $L^2$-norm and $H^1$-norm of the velocity field, respectively. For the $L^2$-norm of the pressure we observe the optimal rate of $k$. It is notable that the convergence behavior of the stabilized formulation considered here is highly regular on all considered meshes, as opposed to the convergence behavior of the non-stabilized FCM formulation; cf. [21].

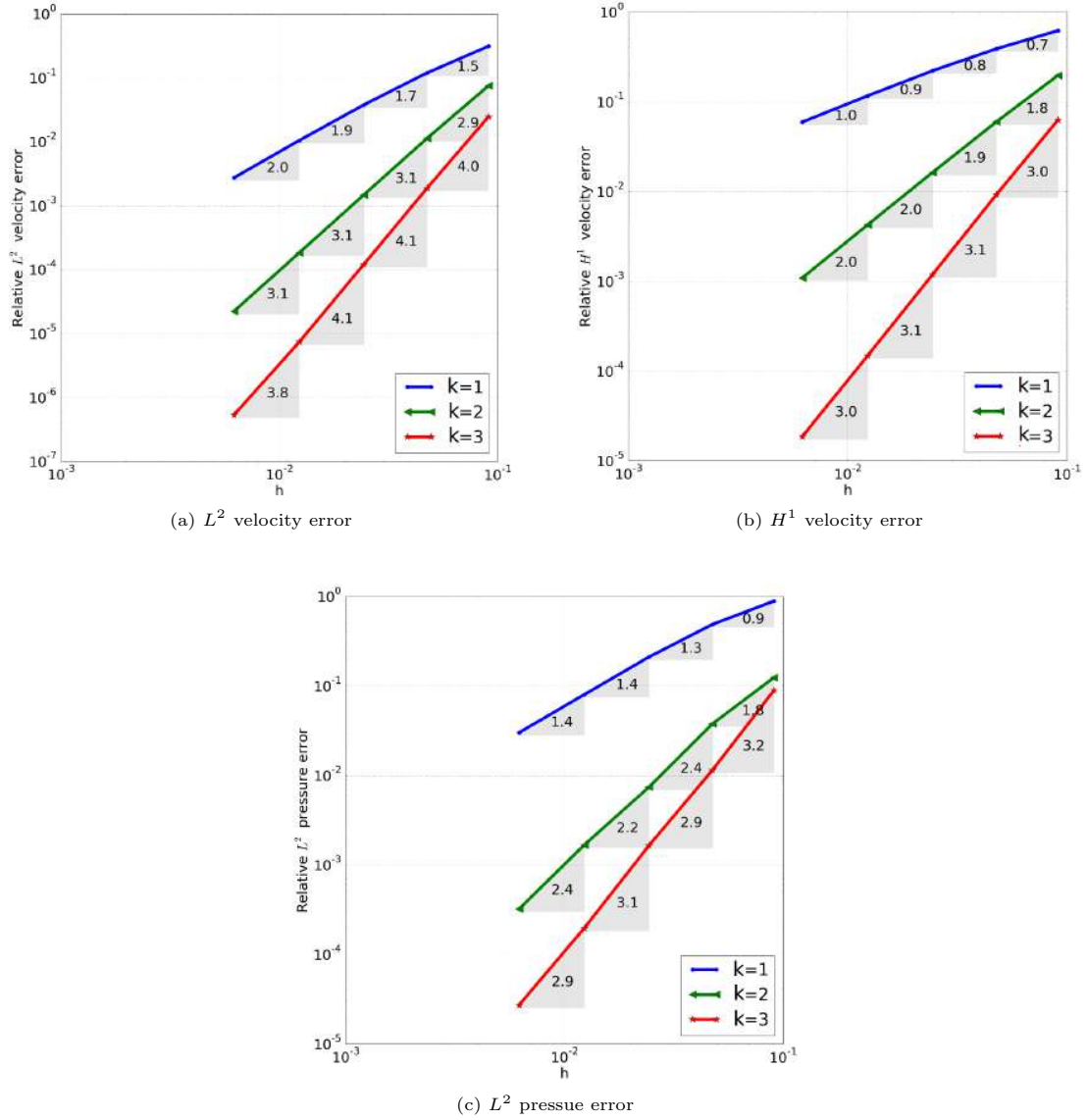

Figure 9: Mesh convergence study of the steady Navier-Stokes equations with $\mu = 1$ in a quarter annulus ring using the skeleton-stabilized formulation with linear, quadratic and cubic full-regularity B-splines.

In Figure 10 we study the solution sensitivity with respect to the skeleton stabilization parameter $\gamma$, where the ghost-penalty parameter is fixed at $\widetilde{\gamma} = 5 \times 10^{-3}$. The $h$-convergence behavior of the solution using $C^1$-continuous quadratic B-splines is studied for a wide range of stabilization parameters, viz. $\gamma \in (5 \times 10^{-5}, 10)$. We observe that the pressure stabilization parameter $\gamma$ does not affect the accuracy of the velocity field in the $L^2$-norm and $H^1$-norm. This behavior is expected, as the skeleton-penalty term acts only on the pressure field, similar as in the conforming isogeometric analysis setting considered in Ref. [22]. The pressure solution accuracy is affected by the selection of the stabilization parameter, but Figure 10 conveys that the parameter can be selected from a wide range (approximately $\gamma \in (5 \times 10^{-4}, 5 \times 10^{-1})$) with minor influence on the accuracy. The convergence rate remains optimal for all considered choices of the stabilization parameter.

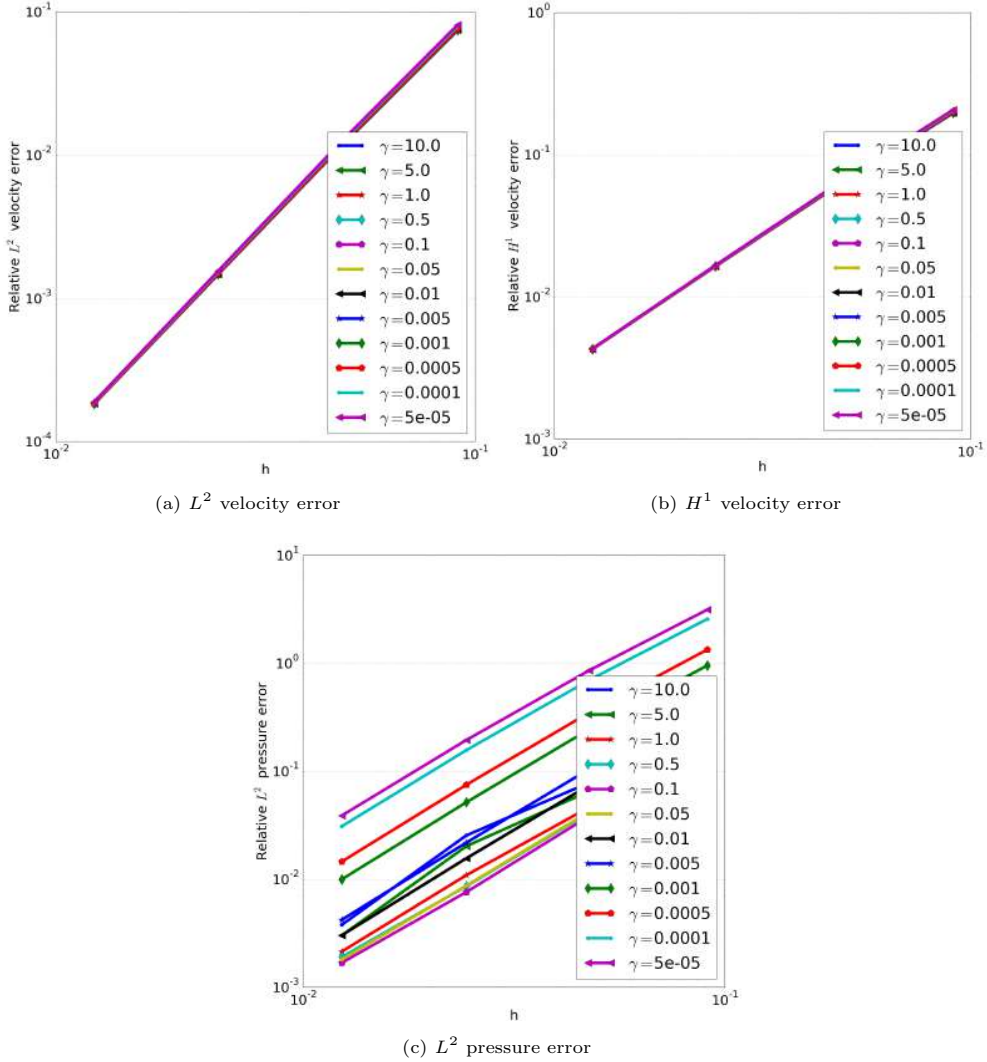

(a) $L^2$ velocity error

(b) $H^1$ velocity error

(c) $L^2$ pressure error

Figure 10: Sensitivity of the quadratic spline approximation of the Navier-Stokes problem with $\mu = 1$ on the quarter annulus ring with respect to the stabilization parameter $\gamma$.

To assess the inf-sup stability of the proposed method we perform a numerical study of the generalized inf-sup constant (see also, *e.g.*, Ref. [39]) similar to that presented for the conforming setting in Ref. [22]. The discrete stability constant associated with the stabilized Stokes system matrix, $\lambda_h$, can be computed as the square root of the smallest non-zero eigenvalue of the generalized eigenvalue problem

$$(\mathbf{B}\mathbf{A}^{-1}\mathbf{B}^{\mathbf{T}} + \mathbf{S})\mathbf{q} = \lambda_h{}^2 \mathbf{M}_{pp}\mathbf{q}, \qquad (23)$$

where $\mathbf{A}$, $\mathbf{B}$, and $\mathbf{S}$ are defined as in Sec. 4, $\mathbf{M}_{pp}$ is the Gramian matrix associated with the pressure basis, *i.e.*, $(\mathbf{M}_{pp})_{ij} = (N_i^p, N_j^p)_{Q^h}$. The discrete norm in the pressure space associated with the Gramian matrix is defined as

$$\left\|q^h\right\|_{Q^h}^2 := \left\|q^h\right\|_{L^2(\Omega)}^2 + \gamma \sum_{F \in \mathcal{F}_{skeleton}^h} \int_F \mu^{-1} h_F^{2k+1} \left|[\![\partial_n^k q^h]\!]\right|^2 \, d\Gamma. \qquad (24)$$

It should be noted that since the norm $\|\cdot\|_{Q^h}$ is stronger than $\|\cdot\|_{L^2(\Omega)}$, numerical inf-sup stability in $\|\cdot\|_{Q^h}$ implies stability for the case that the Gramian matrix $\mathbf{M}_{pp}$ is defined as the $L^2$ pressure mass matrix.

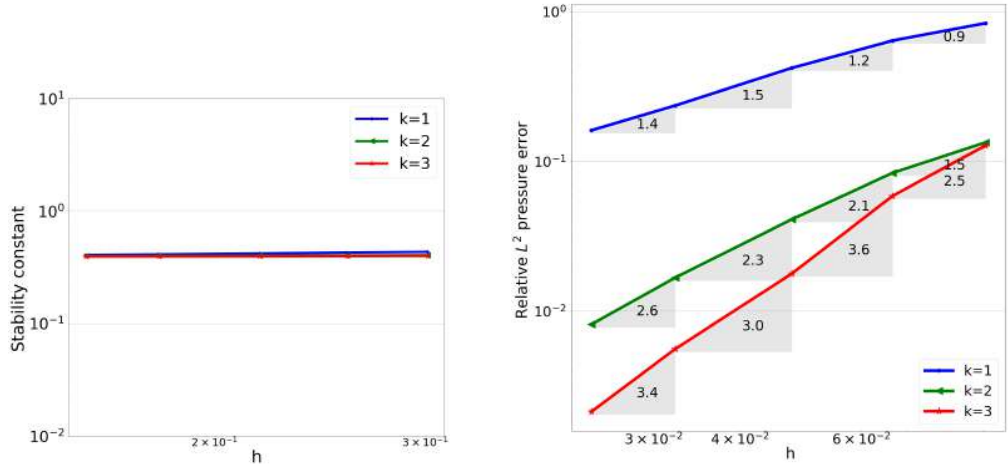

(a) $\gamma = 10$ $(k=1)$, $0.1$ $(k=2)$, $5 \cdot 10^{-4}$ $(k=3)$

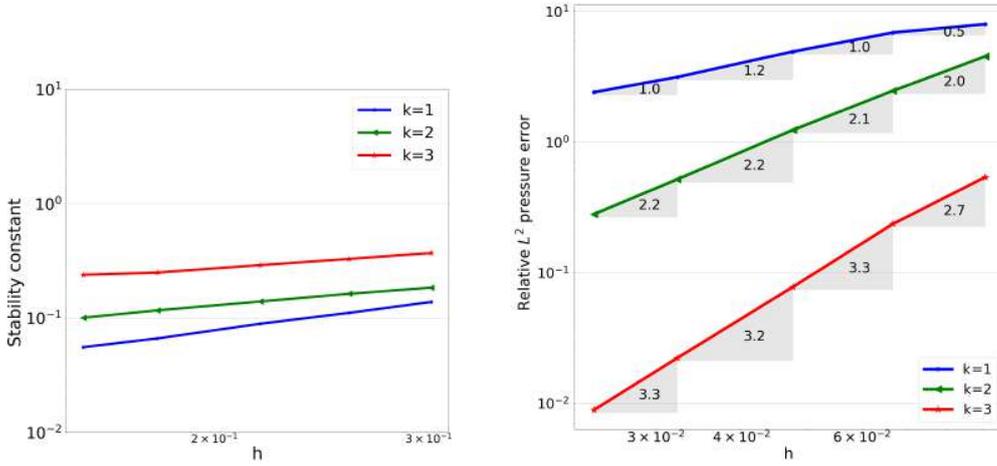

(b) $\gamma = 10^{-5}$ (too small)

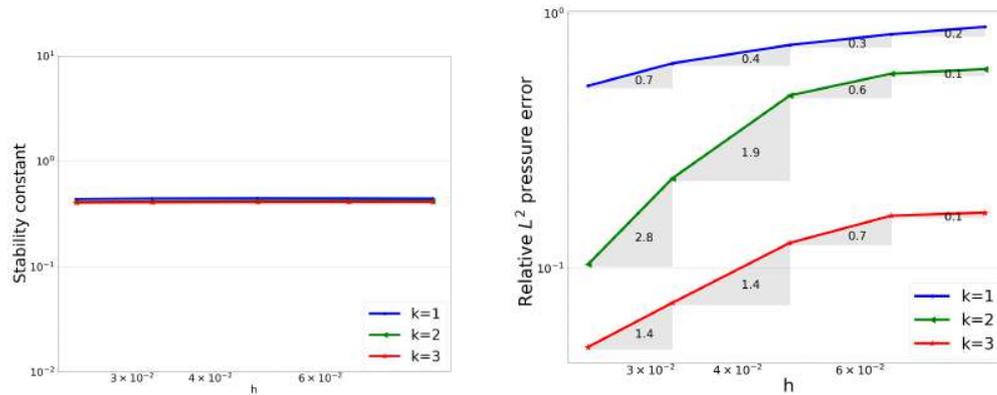

(c) $\gamma = 10^3$ (too large)

Figure 11: Effect of the skeleton-penalty parameter choice $\gamma$: Discrete inf-sup constant (left column) and the $L^2$ pressure error (right column) for various mesh sizes ($h$) and spline degrees ($k = 1, 2, 3$). Figure 11a corresponds to the skeleton-penalty parameter $\gamma = 10$ $(k = 1)$, $0.1$ $(k = 2)$, $5 \times 10^{-4}$ $(k = 3)$, and ghost-penalty parameter $\widetilde{\gamma} = 10^{-k-1}$. Figure 11b considers the case for which the skeleton-penalty parameter is chosen too small, whereas Figure 11c represents a too large selection of the skeleton-penalty parameter.

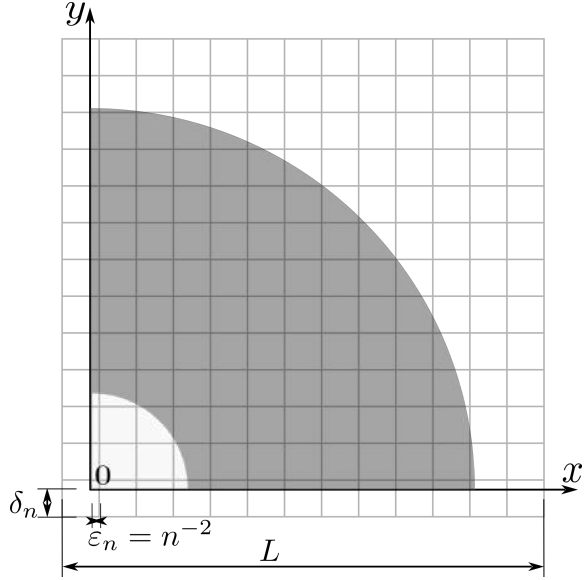

Figure 12: An immersed quarter annulus domain with sliver elements (along the $x$ and $y$ axes) created by offsetting the ambient domain mesh. The square ambient domain of size $L \times L$ is discretized using $n \times n$ elements. The offset $\delta_n$ (in both directions) is selected such that slivers of size $\varepsilon_n = \frac{1}{n^2}$ are created. Note that by this choice the sliver size tends to zero as the number of elements increases.

The effect of the skeleton-penalty parameter $\gamma$ is investigated by computation of the inf-sup constants. The results of this study are shown in Figure 11 for various selections of the skeleton-penalty parameter. The ghost-penalty parameter is fixed at $\widetilde{\gamma} = 10^{-k-1}$ ($k = 1, 2, 3$). When the skeleton-penalty parameter $\gamma$ is chosen as $\gamma = 10$ ($k = 1$), $0.1$ ($k = 2$), $5 \times 10^{-4}$ ($k = 3$), the numerical inf-sup constants are bounded away from 0, independent of the mesh size. Accurate solutions with optimal convergence rates are then obtained. When the skeleton-penalty parameter $\gamma$ is chosen too small (Figure 11b) the discrete inf-sup is affected. The relative pressure errors are observed to increase, although the rates remain optimal. When $\gamma$ is too large (Figure 11c) the system is stable, but the accuracy of the computed results is affected significantly. This effect is most notable for the highest-order case ($k = 3$), which generally requires the lowest value for $\gamma$. As can be observed, the pressure errors are considerably higher than in the case of $\gamma = 5 \times 10^{-4}$, and the asymptotic convergence rates are only attained for small mesh sizes.

To study the effectivity and robustness of the velocity ghost-penalty stabilization, we consider the problem setup shown in Figure 12. The same quarter annnulus as considered above is immersed in an ambient domain of length $L = 5$. The number of elements per direction, $n$, is varied. The offset of the ambient domain with respect to the coordinate system is chosen such that when refining the mesh (increasing $n$), the width of the cut-bands on the left and the bottom of the domain reduces to 0. For elliptic problems, it is well known that such sliver-cut configurations may lead to an unbounded control over flux variables at the cut boundary zone [26, 34]. We here fix the skeleton-penalty parameter at the optimized values from the previous study, viz. $\gamma = 10$ ($k = 1$), $0.1$ ($k = 2$), $5 \times 10^{-4}$ ($k = 3$). The Nitsche stabilization parameter is $\beta = 6(k + 1)^2$, while the mesh-size $h^{-1}$ in equation (14) is of the regular background mesh. The number of elements per edge is $n = 11, 15, 21, 31, 41$. We consider three different cases for the selection of the ghost-penalty parameter: the quasi-optimal (ad-hoc) choice $\widetilde{\gamma} = 10^{-k-1}$ ($k = 1, 2, 3$), a vanishing ghost-penalty $\widetilde{\gamma} = 0$, and a too large value $\widetilde{\gamma} = 10^4$. From Figure 13 we observe that for all three cases the discrete generalized inf-sup constant is bounded away from 0. It should be recalled that the pressure skeleton-penalty term $s^h_{\text{skeleton}}$ is defined on the $\mathcal{F}^h_{\text{skeleton}}$ structure, which includes $\mathcal{F}^h_{\text{ghost}}$ as a subset. For the optimized case (Figure 13a), we obtain optimal convergence rates for both the velocity and pressure. When the ghost-penalty vanishes (Figure 13b), an adverse effect on the $H^1$ error of the velocity is observed, which is in accordance with the literature on elliptic problems [25, 26, 40]. When $\widetilde{\gamma} = 1 \cdot 10^4$ (Figure 13c), a negative impact on the accuracy of the velocity is observed, which also leads to sub-optimal convergence of the pressure in the $L^2$ norm. For lower orders, convergence is even not apparent at all.

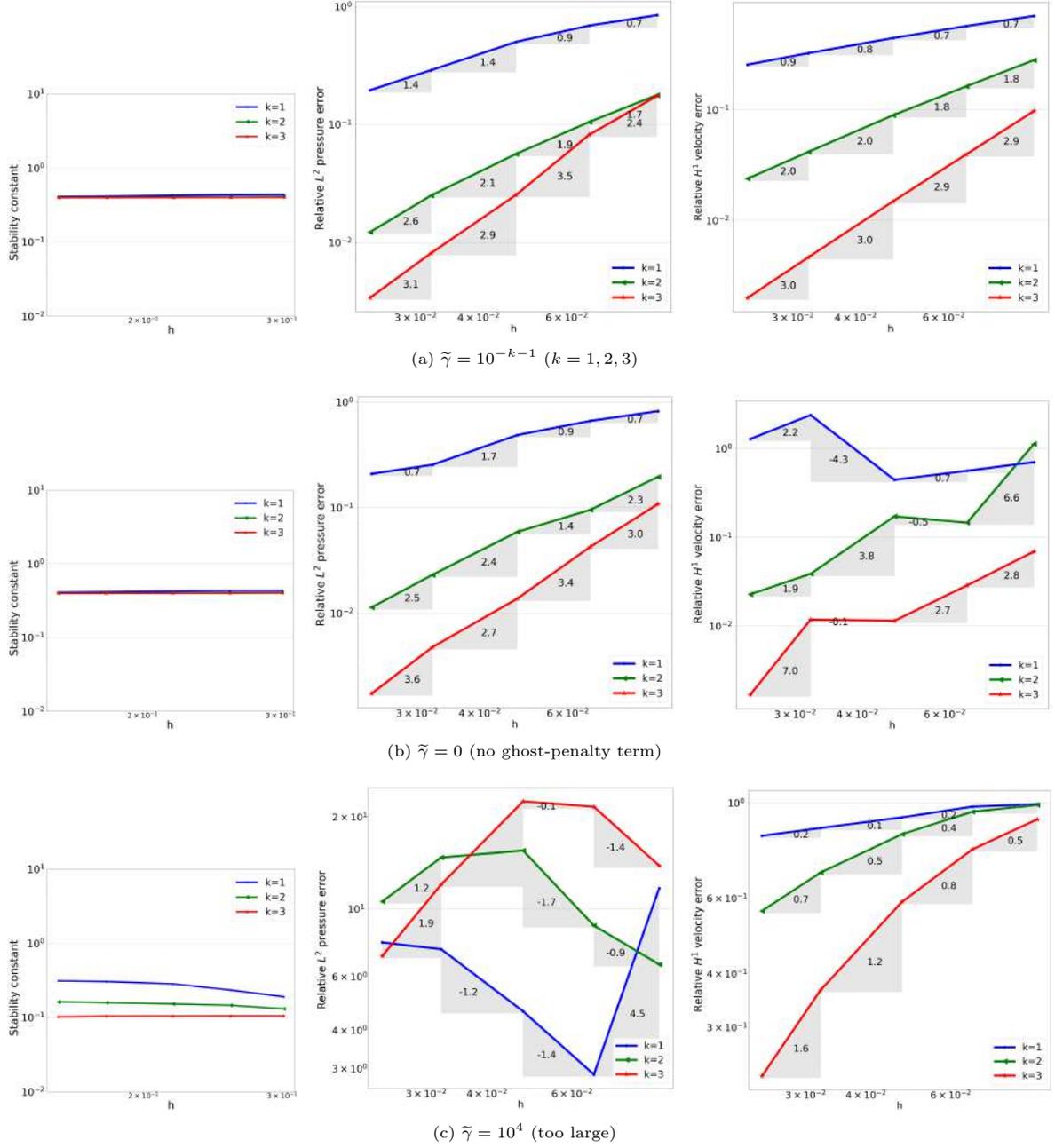

Figure 13: Effect of the ghost-penalty parameter $\widetilde{\gamma}$: Discrete inf-sup constant (first column), the $L^2$ pressure error (middle column) and the $H^1$ velocity error (last column) for various mesh sizes ($h$) and spline degrees ($k = 1, 2, 3$). The skeleton-penalty parameter is fixed at $\gamma = 10$ ($k = 1$), $0.1$ ($k = 2$), $5 \times 10^{-4}$ ($k = 3$). Figure 13a corresponds to the ghost-penalty parameter $\widetilde{\gamma} = 10^{-k-1}$. Figure 13b considers the case for which the skeleton-penalty parameter is 0, whereas Figure 13c represents a too large selection of the ghost-penalty parameter.

## 5.2. Navier-Stokes flow around a cylinder

In this section we revisit the benchmark problem proposed by Schäfer and Turek [41], which we considered in the mesh-conforming isogeometric analysis in Ref. [22]. To perform a detailed assessment of the convergence behavior of the proposed formulation, in Section 5.2.1 we first consider a simplified cylinder flow with a smooth manufactured solution. Subsequently, in Section 5.2.2 we validate our formulation using the benchmark problem.

### 5.2.1. Manufactured solution

We consider steady Navier-stokes flow in a unit square with a cylinder of radius $R = \frac{1}{8}$ located at the center of the square. The viscosity is taken as $\mu = 1$. Velocity field Dirichlet conditions that are compatible with the manufactured solution [3]

$$\mathbf{u} = \begin{pmatrix} 2e^x(-1+x)^2 x^2(y^2-y)(-1+2y) \\ (-e^x(-1+x)x(-2+x(3+x))(-1+y)^2 y^2) \end{pmatrix} \tag{25a}$$

$$p = (-424 + 156e + (y^2-y)(-456 + e^x(456 + x^2(228 - 5(y^2-y)) + 2x(-228 + (y^2-y)) \\ + 2x^3(-36 + (y^2-y)) + x^4(12 + (y^2-y)))))  \tag{25b}$$

are imposed on the complete boundary. A Lagrange multiplier is introduced to enforce the average pressure condition. As quantities of interest we consider the lift and drag coefficients, $c_L$ and $c_D$, respectively. These coefficients are defined as

$$c_D = \frac{F_D}{\rho \bar{U}^2 R} = \frac{\mathcal{R}(\mathbf{u}, p; \boldsymbol{\ell}_1)}{\rho \bar{U}^2 R}, \qquad c_L = \frac{F_L}{\rho \bar{U}^2 R} = \frac{\mathcal{R}(\mathbf{u}, p; \boldsymbol{\ell}_2)}{\rho \bar{U}^2 R}, \tag{26}$$

where $F_D$ and $F_L$ are the resultant drag and lift forces acting on the cylinder, which are evaluated weakly as (see *e.g.*, [42–44])

$$\mathcal{R}(\mathbf{u}, p; \boldsymbol{\ell}_i) := c(\mathbf{u}; \mathbf{u}, \boldsymbol{\ell}_i) + a(\mathbf{u}, \boldsymbol{\ell}_i) - 2\mu \langle \nabla^s \boldsymbol{\ell}_i \cdot \mathbf{n}, \mathbf{u} - \mathbf{g} \rangle_{\Gamma_D} + b(p, \boldsymbol{\ell}_i) - (\mathbf{f}, \boldsymbol{\ell}_i),$$

with $\boldsymbol{\ell}_i \in [H^1_{0,\partial\Omega\setminus\Gamma}(\Omega)]^d$ and $\boldsymbol{\ell}_i|_\Gamma = -\mathbf{e}_i$, $i = 1, 2$. For the manufactured solution considered here we select $\rho\bar{U}^2 = 1$, where $\rho = 1$ is the fluid density, and $\bar{U} = 1$ is a representative velocity scale.

To study the influence of the geometry representation on the convergence behavior of the quantities of interest (26) we consider quadratic B-spline discretizations with two levels of bi-sectioning, *viz.* $\varrho_{\max} = 1$ and $\varrho_{\max} = 7$. The stabilization parameters are chosen as $\gamma = \widetilde{\gamma} = 10^{-3}$. Since the bi-sectioning level is defined relative to the untrimmed cells, the geometry is refined along with the ambient domain mesh refinement. In Figure 14 the computed solutions for the two settings of the maximum bi-sectioning level are shown on three different meshes with $h = \frac{1}{3}, \frac{1}{5}, \frac{1}{9}$. The top panels pertain to the case of only a single level of bi-sectioning, *i.e.* $\varrho_{\max} = 1$, for which the geometric mismatch of the domain is clearly visible on all considered meshes. The observed variations of the solution under mesh refinement can be attributed to both the geometric error and discretization error being resolved simultaneously. In contrast, for the case of $\varrho_{\max} = 7$ as shown in the bottom panels, the geometry mismatch is not discernible on any of the meshes. In this case, the error in the approximation is dominated by the discretization.

In Figure 15 we study the convergence behavior of the errors in the lift and drag coefficients for the two settings of the maximum bi-sectioning level. For the case of $\varrho_{\max} = 1$ we observe both quantities of interest to converge with an approximate rate of 2, whereas a rate of approximately 4 is observed in the case of $\varrho_{\max} = 7$. The reduced rate of convergence for $\varrho_{\max} = 1$ (in comparison to the $\varrho_{\max} = 7$ case) can be attributed to the fact that the geometry is only first order accurate, *i.e.*, on the maximum bi-sectioning level the trimmed boundary is interpolated linearly. The observed rate of 2 indicates that the error for $\varrho_{\max} = 1$ is in this case dictated by the geometry representation. When the geometry is properly resolved (in this case by selecting $\varrho_{\max} = 7$) the observed convergence rates are consistent with the double-order ($2k = 4$) super-convergence behavior for smooth functionals that is generally observed for conforming Galerkin approximations [22, 45]. Under further mesh refinement the asymptotic convergence rate will inevitably reduce to the geometry-dictated rate of 2. For the results presented here, on the finest mesh considered the geometry error for $\varrho_{\max} = 7$ is still negligible compared to the discretization error.

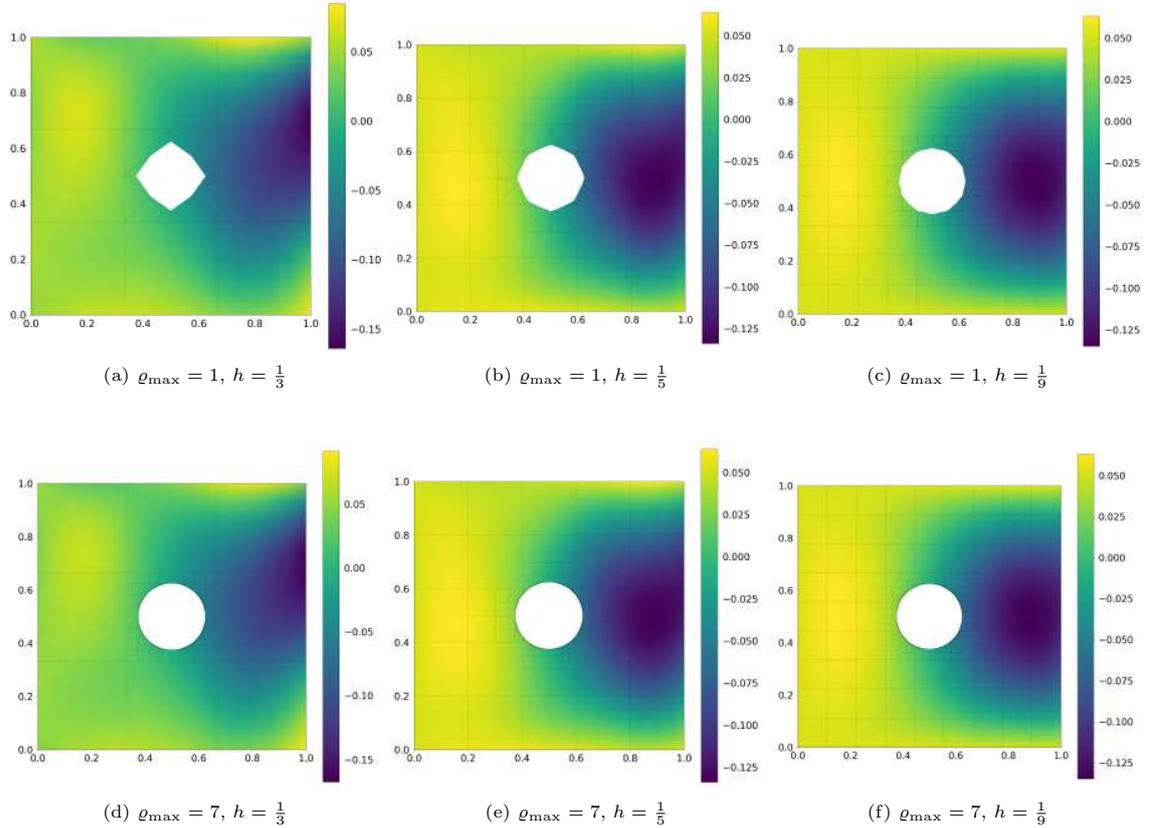

Figure 14: Solutions for the cylinder flow problem with the manufactured solution (25) on a sequence of meshes with (top row) one level of bi-sectioning ($\varrho_{\max} = 1$) and (bottom row) seven levels of bi-sectioning ($\varrho_{\max} = 7$).

*5.2.2. CFD benchmark*

We now consider the benchmark problem proposed by Schäfer and Turek [41]. In this test case a cylinder of radius $R = 0.05 \, \mathrm{m}$ is placed in a channel of height $H = 0.41 \, \mathrm{m}$ and length $L = 2.2 \, \mathrm{m}$. The center of the cylinder is positioned at a horizontal distance of $W = 0.2 \, \mathrm{m}$ from the inflow boundary at $x = 0$, and at a vertical distance of $W = 0.2 \, \mathrm{m}$ from the bottom channel wall at $y = 0$. Note that the cylinder has a small offset with respect to the center line of the channel, introducing an asymmetry in the problem. At the inflow boundary the parabolic flow profile

$$\mathbf{u}(0, y) = \begin{pmatrix} 4U_m y(H - y)/H^2 \\ 0 \end{pmatrix}$$

with maximum velocity $U_m$ is imposed. There is no slip at the bottom and top boundaries, as well as along the surface of the cylinder, and a natural boundary condition is used at the outflow boundary ($x = L$). The kinematic viscosity of the fluid is taken as $\mu = 1 \times 10^{-3} \, \mathrm{m}^2/\mathrm{s}$.

We consider the case of Re = 20 – with the Reynolds number defined as Re $= 2\bar{U}R/\mu$ (with mean inflow velocity $\bar{U} = \frac{2}{3}U_m$) – for which a steady flow is obtained. We consider a sequence of uniform refinements of a relatively coarse ambient domain mesh consisting of $36 \times 22$ elements. This coarsest mesh is a non-uniformly spaced full-regularity B-spline patch, with the knot values selected so that relatively small elements are obtained in the neighborhood of the cylindrical inclusion (see Figure 16). The outer boundaries of the ambient domain mesh coincide with the boundaries of the physical domain. The essential boundary conditions are, however, still enforced weakly by Nitsche's method. Figure 16 shows the speed and pressure computed on the three times refined second-order B-spline mesh, which results in a system of $n_{\text{dof}} = 148476$ degrees of freedom. The obtained result is visually in good agreement with the benchmark result, and is free of pressure oscillations near the immersed boundary of the cylinder.

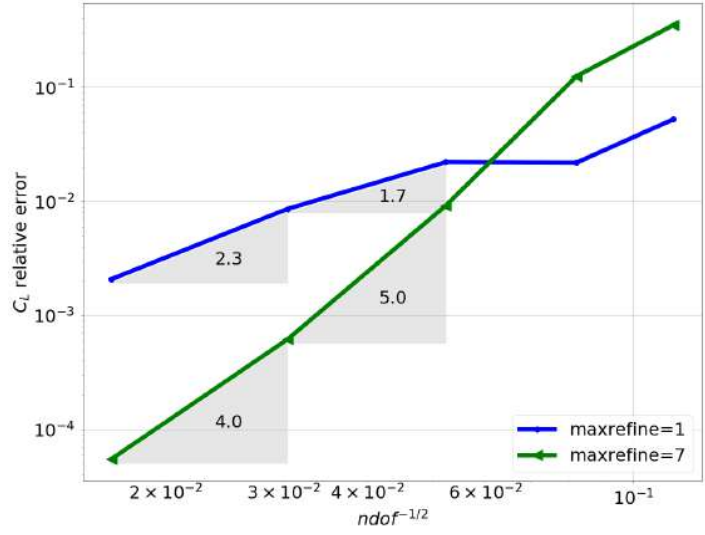

(a) Lift coefficient, $c_L$

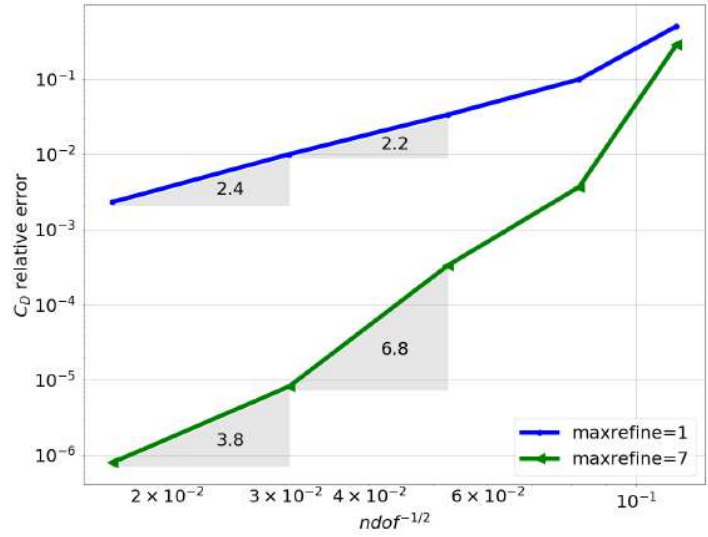

(b) Drag coefficient, $c_D$

Figure 15: Convergence behavior of the lift and drag coefficients for the cylinder flow problem with the manufactured solution (25) with two settings for the maximum bi-sectioning level.

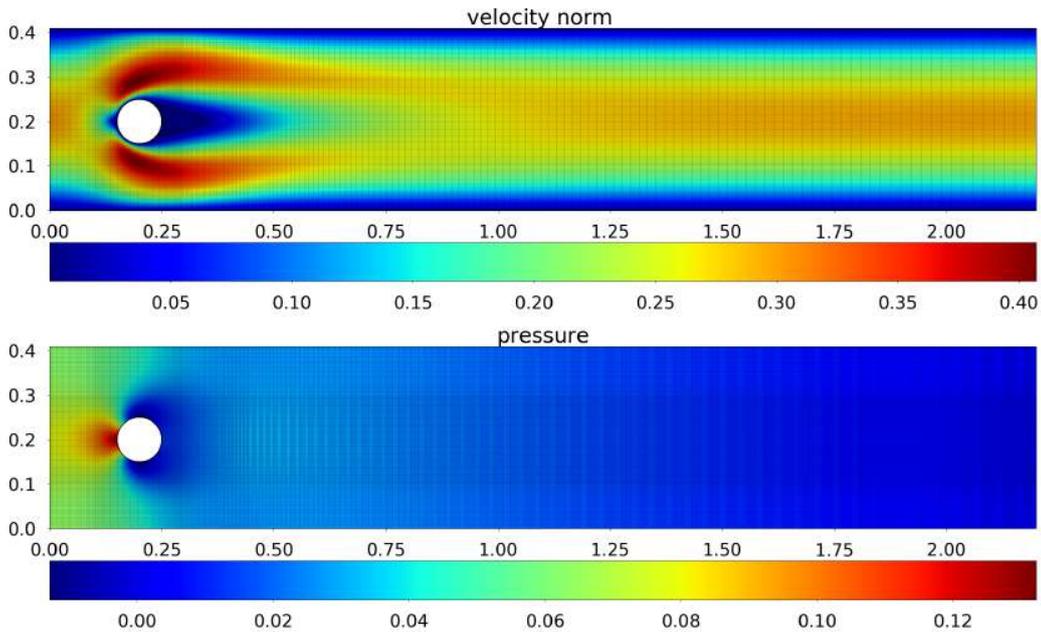

Figure 16: Velocity magnitude (top) and pressure (bottom) solutions of the steady cylinder flow problem using quadratic B-splines with $n_{\text{dof}} = 148476$.

In Table 1 we present the mesh converge results for various quantities of interest, *viz.* the lift and drag coefficients, $c_L$ and $c_D$, respectively, and the pressure drop over the cylinder, $\Delta p$. The drag and lift coefficients are defined as in equation (26), and the pressure drop is defined as $\Delta p = p(W - R, W) - p(W + R, W)$.

From Table 1 it is observed that on the finest mesh all quantities of interest are in excellent agreement with the benchmark result [46]. We note that – despite the fact that we here have selected the bisectioning integration depth to ten – with the higher-order approximation of these quantities of interest we anticipate deterioration of the approximation properties associated with the reduced geometric regularity of the immersed boundary approximation. We expect that the loss of convergence rate especially observed for the drag coefficient can be attributed to this, but further investigation of this aspect is warranted.

| Level | $n_{\text{dof}}$ | $h_{\text{cylinder}}$ | $C_D$ | $C_L$ | $\Delta p$ |
|---|---|---|---|---|---|
| 0 | 2724 | 0.0142857 | 5.85317738567 | 0.009231742514 | 0.18749814772 |
| 1 | 9984 | 0.0071428 | 5.57503959343 | 0.011002907897 | 0.12424885895 |
| 2 | 38052 | 0.0035714 | 5.57961186549 | 0.010542479971 | 0.11504949327 |
| 3 | 148476 | 0.0017857 | 5.57989543774 | 0.010575002816 | 0.11703150638 |
| Ref. [46] | | | 5.57953523384 | 0.010618948146 | 0.11752016697 |

Table 1: Computed values of the drag and lift coefficients and pressure drop on four refinement levels using the skeleton-stabilized formulation with quadratic B-splines. The mesh resolution in the proximity of the cylinder is indicated by $h_{\text{cylinder}}$, and the total number of degrees of freedom by $n_{\text{dof}}$.

### 5.3. Scan-based analysis of a porous medium flow

To demonstrate the potential of the proposed formulation for geometrically and topologically complex three-dimensional domains, we consider a creeping flow through a porous medium. The porous medium under consideration is made of sintered glass beads. Three-dimensional gray-scale voxel data of the specimen is obtained by a $\mu$CT-scanner with a voxel resolution of $25\,\mu$m. We here consider a representative domain of $50 \times 50 \times 50$ voxels. The size of this cubic domain is denoted by $L = 1.25$ mm.

In this numerical simulation we compare the immersogeometric approach considered in this work with a voxel-based analysis, which is commonly the method of choice for this type of analyses. The geometric model for the voxel-based analysis is obtained by direct segmentation of the gray-scale data (Figure 17a), where all gray-scale values larger than a calibrated threshold are eliminated from the domain. This voxel data threshold is calibrated so that the porosity of the segmented void space is as close as possible to the specified target porosity of 28%. The corresponding voxel model for the void space is shown in Figure 17b. It is noted that as a consequence of the discrete character of the gray-scale data the realized void-space porosity is equal to 28.1%.

B-spline based finite cell domains are obtained using the procedure proposed in Ref. [19]. First the gray-scale voxel data is smoothened by convolution of that data with a (second-order) B-spline basis constructed over the voxel grid, an operation that bears resemblance with Gaussian blurring. Then a relatively coarse (second-order) B-spline mesh is created over the ambient domain matching the scan size, so that the outer boundaries of the pore space reside in the boundaries of the scan domain. The smooth B-spline level set function is then segmented using the bi-sectioning procedure described in Section 2.2, where the threshold is calibrated based on the required porosity of 28%. Note that, in contrast to the voxel data, the target porosity can be matched up to a tolerance specified by the user (in this case up to six digits). Figures 17c and 17d show the geometric models obtained using maximum bi-sectioning depths of $\varrho_{\text{max}} = 1, 3$ with respect to an ambient domain mesh consisting of $12 \times 12 \times 12$ elements. Evidently, the geometry representation becomes smoother as the bi-sectioning depth is increased, indicating that the underlying smoothened level set function is better resolved by the integration mesh. The improvement of the geometry representation can also be observed from Figure 18a, which conveys that the threshold correction following from the porosity-calibration procedure decreases as the bi-sectioning depth increases.

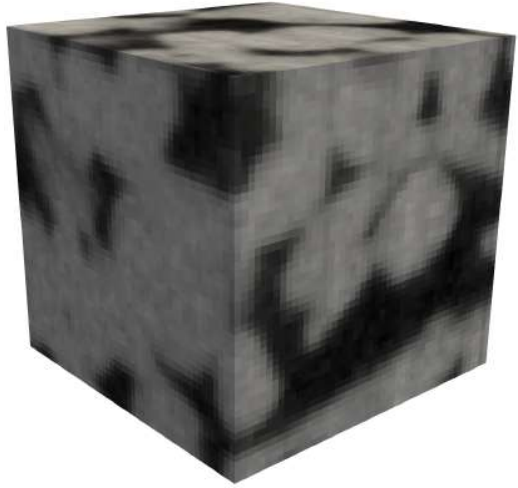
(a) Gray-scale scan data

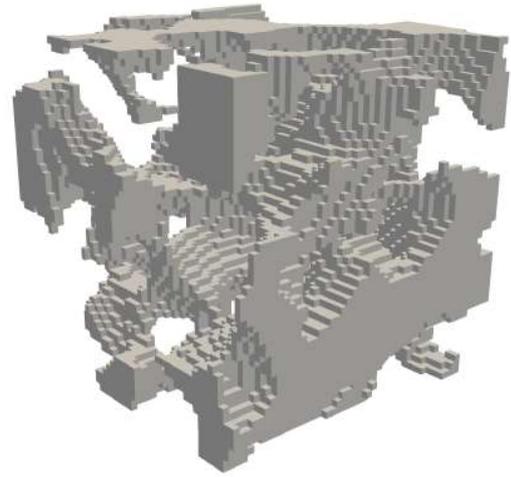
(b) Voxel domain

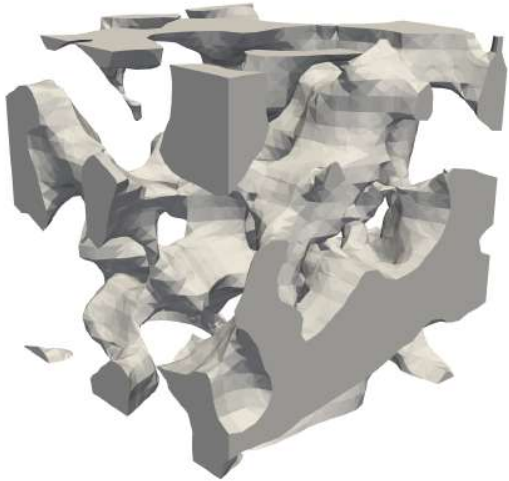
(c) Immersogeometric: $\varrho_{\max} = 1$

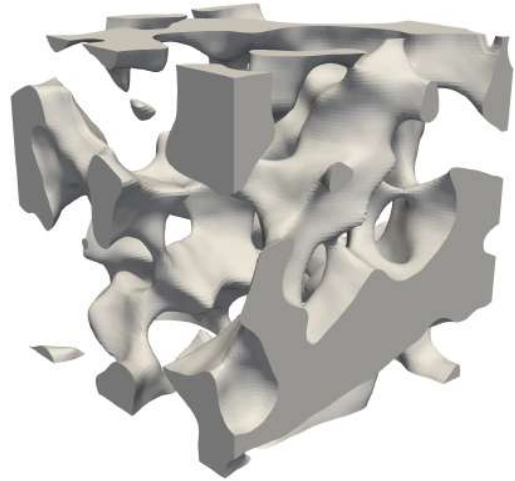
(d) Immersogeometric: $\varrho_{\max} = 3$

Figure 17: Various $\mu$CT-based geometric models of a sintered glass beads porous medium specimen. The original scan data consists of 125000 gray-scale voxels with a resolution of $25\,\mu$m. The illustrated immersed domains correspond to a $12 \times 12 \times 12$ background mesh with different maximum bi-sectioning levels. The porosity of all shown domains is calibrated to 28%.

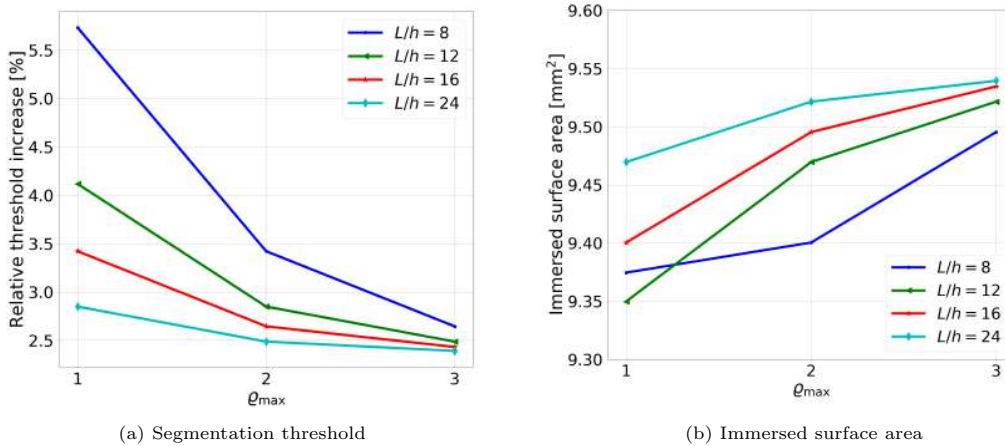

(a) Segmentation threshold  (b) Immersed surface area

Figure 18: Porosity-calibrated gray-scale threshold offset and corresponding immersed surface area for the sintered glass beads specimen with various background meshes and bi-sectioning depths. The threshold increase is defined relative to the calibrated voxel data threshold.

Comparison of the voxel model and the immersogeometric model reveals that both geometric models are visually very similar in terms of micro-structural features, which is expected based on the calibrated porosity. Evidently, the surface representation of the models is completely different. Whereas a staircase representation of the original gray-scale scan data is obtained in the voxel model, a multi-resolution segmentation based on the smoothened gray-scale data is obtained in the immersogeometric case. In this multi-resolution representation, the intra-element curvature of the surface is recovered by the bi-sectioning operation (*i.e.*, the geometric linearization error is associated with the size of the integration sub-cells, and not with that of the size of the computational background mesh). This difference in surface representation translates directly in a significant difference of the internal surface area, which is equal to $15.5\,\text{mm}^2$ for the voxel model, and to approximately $9.5\,\text{mm}^2$ for the immersogeometric model. Indeed, a significantly higher surface area is expected in the voxel representation. From Figure 18b it is observed that the immersed surface area can be captured accurately by the bi-sectioning procedure independent of the background mesh size. Adequate representation of the surface is critical in many situations, for example when surface reactions are considered such as in the case of biofilm growth and mineral dissolution/precipitation in porous media [47], or when one is interested in contact line dynamics for multi-phase porous media flows or elasto-capillarity [48–50].

We consider the simulation of a creeping flow governed by the Stokes equations with viscosity $\mu = 10^{-3}\,\text{Pa}\cdot\text{s}$. The flow through the porous medium is forced by imposition of a normal traction difference of $1\,\text{Pa}$ between the inflow (left) and outflow (right) boundaries. All velocity components are zero on the other lateral boundaries, and no slip conditions are imposed on the interior surface. As a reference we consider the voxel method results shown in Figure 19, which have been computed using FLUENT. The settings of this finite volume method – where the degrees of freedom are closely related with the 35,383 pore space voxels – have been tailored to the problem class under consideration. As a quantity of interest we consider the effective permeability, $\kappa = \frac{\mu Q}{L \Delta p}$, with $Q$ the fluid flow discharge. The pressure drop $\Delta p$ is computed as the difference between the average pressure over the inflow domain and that over the outflow domain, which, for all simulations considered in this section is within a few percent of the imposed normal traction difference of $1\,\text{Pa}$. The effective permeability corresponding to the voxel method results in Figure 19 is $\kappa = 8.9 \cdot 10^{-5}\,\text{mm}^2$.

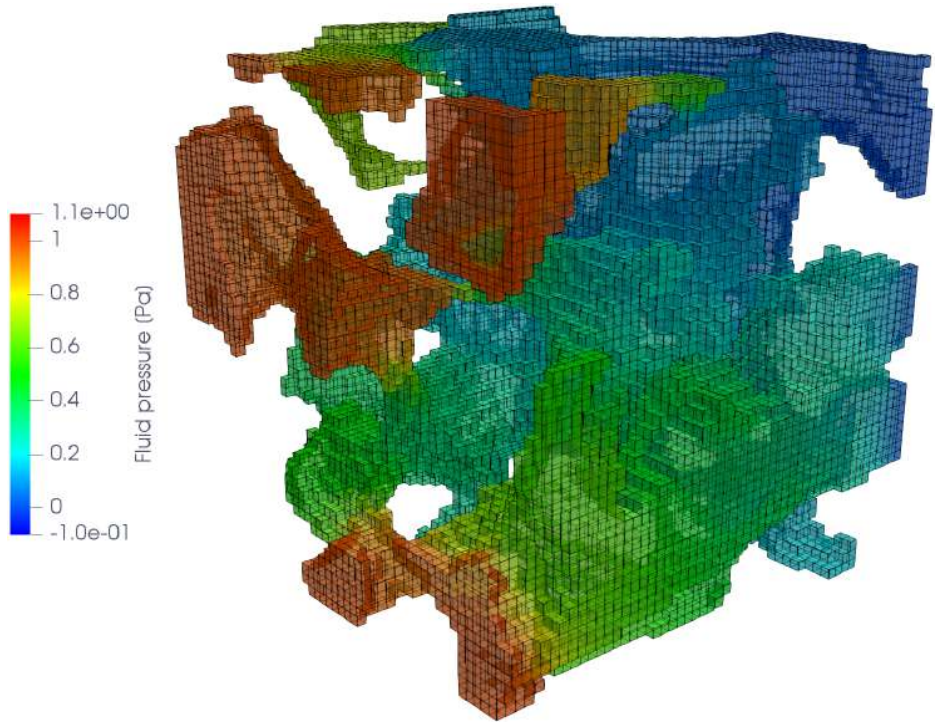

(a) Pressure field

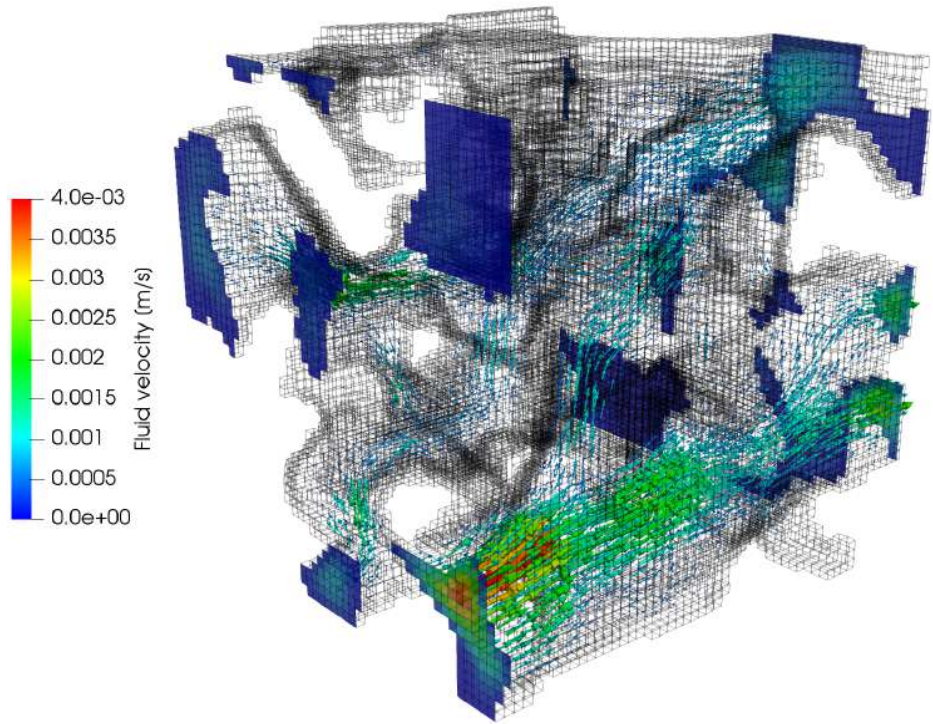

(b) Velocity field

Figure 19: Pressure field and velocity field obtained using the voxel method (computed using FLUENT) for the sintered glass beads specimen.

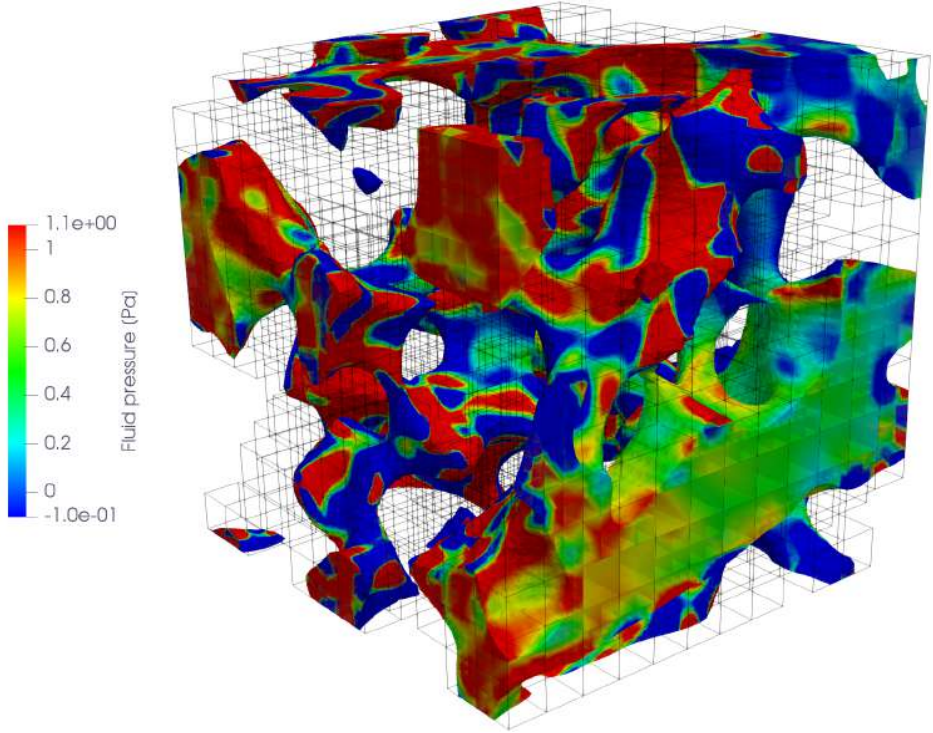

Figure 20: Pressure field computed using immersogeometric analysis with the $C^1$ Taylor-Hood element (cubic velocity, quadratic pressure) without the skeleton-penalty term $s^h_{\text{skeleton}}$ (the ghost-penalty term $s^h_{\text{ghost}}$ is kept for the stability of the Nitsche boundary condition imposition). The unphysical pressure oscillations propagate globally, this in contrast to the result for the topologically simple case in Figure 1. This behavior is a direct consequence of the fact that in this case virtually all elements are cut (i.e., no interior-cells, $\mathcal{F}^h_{\text{ghost}} \simeq \mathcal{F}^h_{\text{skeleton}}$). Here the mesh is $12 \times 12 \times 12$ elements with maximum bi-sectioning depth $\varrho_{\max} = 2$.

Before we discuss the results obtained using the skeleton-stabilized approach proposed in this manuscript, we would like to stress that standard inf-sup stable space approaches without additional stabilization are unsuitable for complex geometries such as the porous medium considered here. This is elucidated in Figure 20, which shows the pressure field computed with the isogeometric $C^1$ Taylor-Hood element using a cubic velocity approximation and a quadratic pressure approximation. The observed global unphysical spurious pressure result is a direct consequence of the fact that for a domain of this kind virtually all elements are cut, which renders direct application of cell-based approaches such as inf-sup stable spaces and GLS-type stabilization ineffective. Refining the mesh in order to resolve the unphysical oscillations will lead to prohibitively large problem sizes, and hence conflicts with the general aim of immersed methods (viz., reducing computational effort when working with complex geometries).

Figure 21 shows the skeleton-stabilized immersogeometric result obtained on an ambient domain mesh of $12 \times 12 \times 12$ elements (with a maximum bi-sectioning depth of $\varrho_{\max} = 3$) using second-order ($k = 2$) B-spline spaces. This immersogeometric model contains approximately 13 times fewer degrees of freedom than the FLUENT simulation in Figure 19. The Nitsche stabilization parameter is taken as $\beta = 100$, and the skeleton-penalty and ghost-penalty parameters are taken as $\gamma = 5 \cdot 10^{-2}$ and $\widetilde{\gamma} = 5 \cdot 10^{-3}$, respectively. Comparison of the immersogeometric model and voxel model shows that both models are in good correspondence, despite the difference in computational resolution. Let us note that the computed flow fields are very similar, although the illustrated flow patterns reveal some visual differences (also because of post-processing artefacts).

In Figure 22 the dependence of the effective permeability on the resolution of the background mesh and the maximum bi-sectioning depth is studied. It is observed that the effect of the bi-sectioning depth decreases as the mesh is refined, which is a consequence of the fact that the geometry bi-sectioning is defined relative to the background mesh. While $\varrho_{\max} = 1$ yields a relatively poor approximation of the geometry on the coarsest ($8 \times 8 \times 8$) mesh, it yields a much better geometry approximation for the finest ($24 \times 24 \times 24$) mesh. Under mesh refinement, the effective permeability is, in general, observed to approach the permeability of the voxel method from above. This observed overestimation of the permeability using the immersogeometric analysis on coarse meshes (specifically the $12 \times 12 \times 12$ mesh) can also be seen from Figure 21b, which, compared to the voxel result in Figure 19b shows a higher flux especially in the upper part of the specimen. It is noted that although the results computed on the coarsest ($8 \times 8 \times 8$) mesh are reasonable, this mesh is too coarse to assign meaning to the results in terms of mesh convergence behavior. The observed behavior of the permeability under mesh refinement indicates that in order to accurately resolve the flow field for the purpose of determining the permeability within a 10% margin, a resolution similar to that of the employed voxel grid is required. This resolution is, however, not required globally, as only a marginal part of the domain contributes significantly to the effective flow through the specimen. In this regard, it is important to note that in terms of computational accuracy the immersogeometric method and the voxel method are fundamentally different. In the case of the voxel method, the computational resolution is intrinsically tied to the geometric model, whereas for the immersogeometric approach, the computational mesh resolution can be controlled independent of the geometric model. As demonstrated in Ref. [19] for an elasticity problem, the decoupling of the computational resolution from the geometric model opens the doors to performing (goal-)adaptive analyses with optimized meshes. In our future work we aim at applying a similar strategy to optimally compute effective permeabilities by only using small elements in regions that contribute to this quantity of interest.

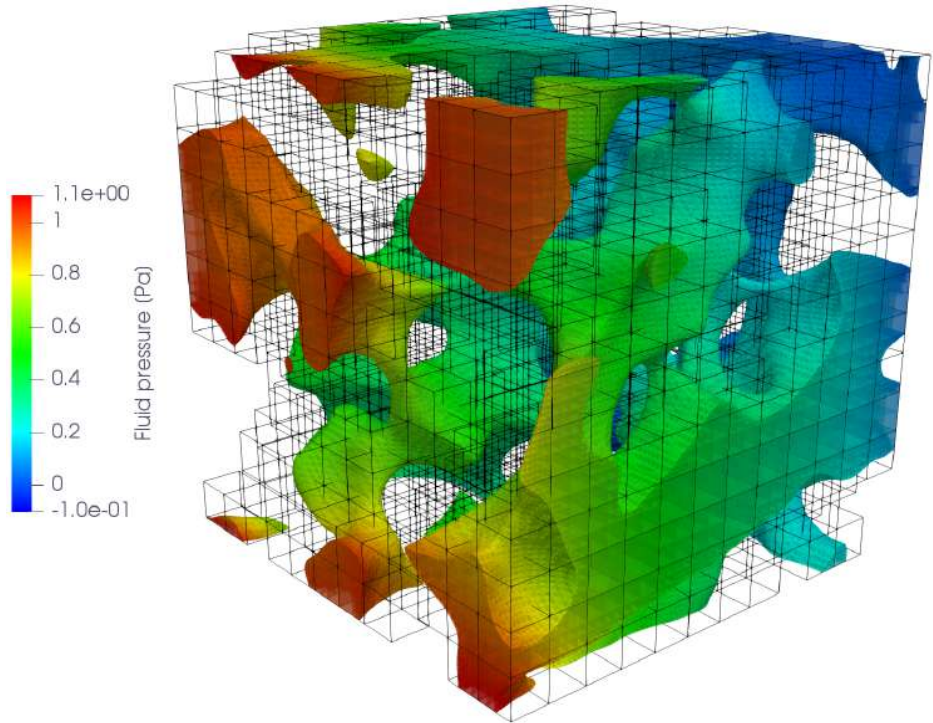

(a) Pressure field

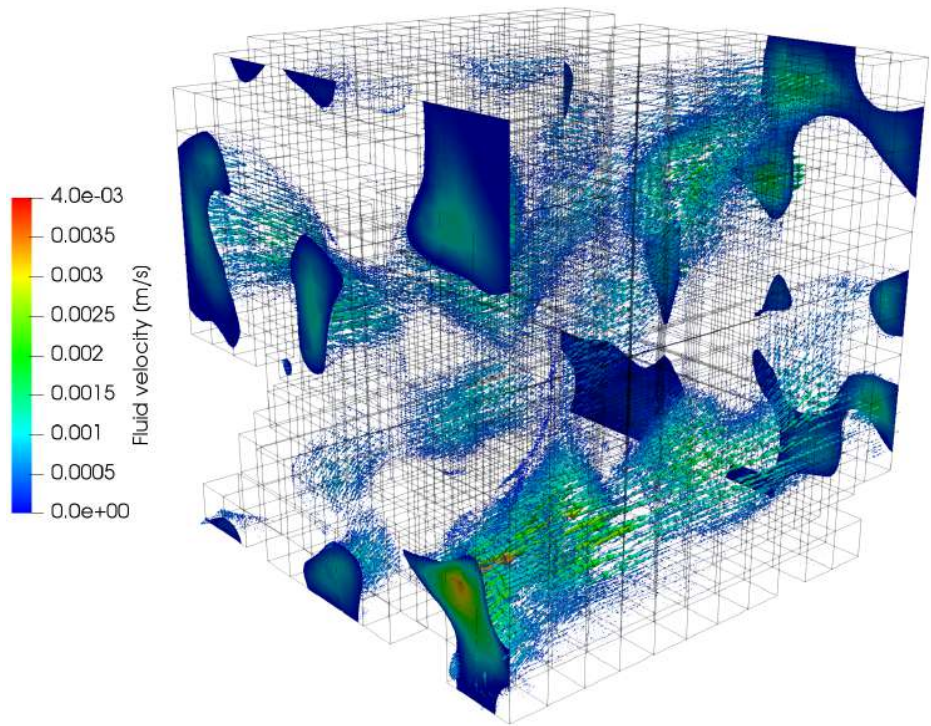

(b) Velocity field

Figure 21: Pressure field and velocity field computed using skeleton-stabilized immersogeometric analysis for the sintered glass beads specimen. The result are based on a $12 \times 12 \times 12$ mesh with maximum bi-sectioning depth $\varrho_{\max} = 3$ and second-order $C^1$ B-splines.

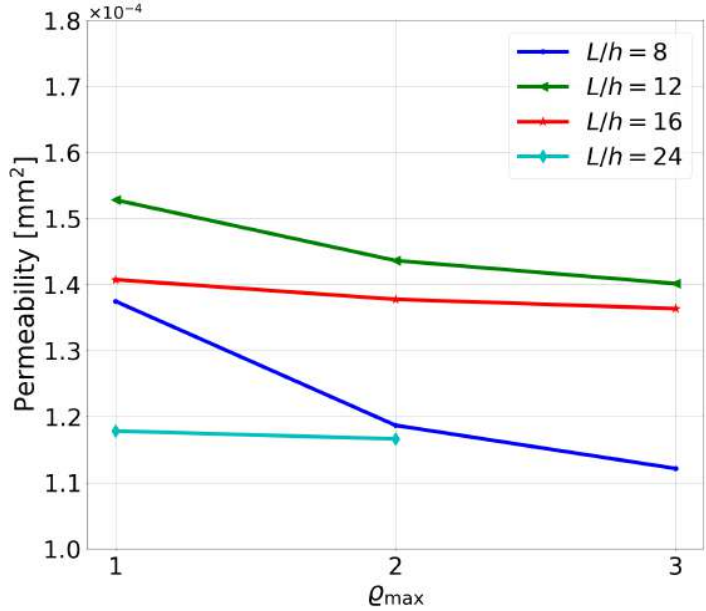

Figure 22: Effective permeability of the sintered glass beads specimen as computed by the skeleton-stabilized immersogeometric analysis using various background meshes and bi-sectioning depths. We note that because of restrictions coming from our current implementation of the (parallel) integration procedure for the trimmed surface integrals associated with the Nitsche terms in the left hand side, the simulation result corresponding to $24 \times 24 \times 24$ elements with maxrefine=3 is absent in this figure.

To demonstrate the potential of the proposed framework for upscaling to bigger simulations, in Figure 23 we show the domain and pressure field result computed for a scan consisting of 2 million voxels. The shown result is computed using $C^1$ continuous quadratic B-splines defined over a $16 \times 16 \times 32$ ambient domain mesh with a maximum bi-sectioning depth of $\varrho_{\max} = 1$. This simulation result illustrates the robustness of the employed geometry parametrization technique. Moreover, the impact of the stabilization terms on the memory consumption of the simulation has been observed to remain limited for this scale of problems, which is a direct consequence of the minimal bandwidth property of the skeleton-stabilized immersogeometric system when using highest regularity B-splines. Further optimization of the employed cut cell integration procedure in terms of computational effort – which, per-element, is significantly smaller for the voxel method – is required to fully unlock the immersogeometric analysis potential at this problem scale and to allow for upscaling to even larger specimens.

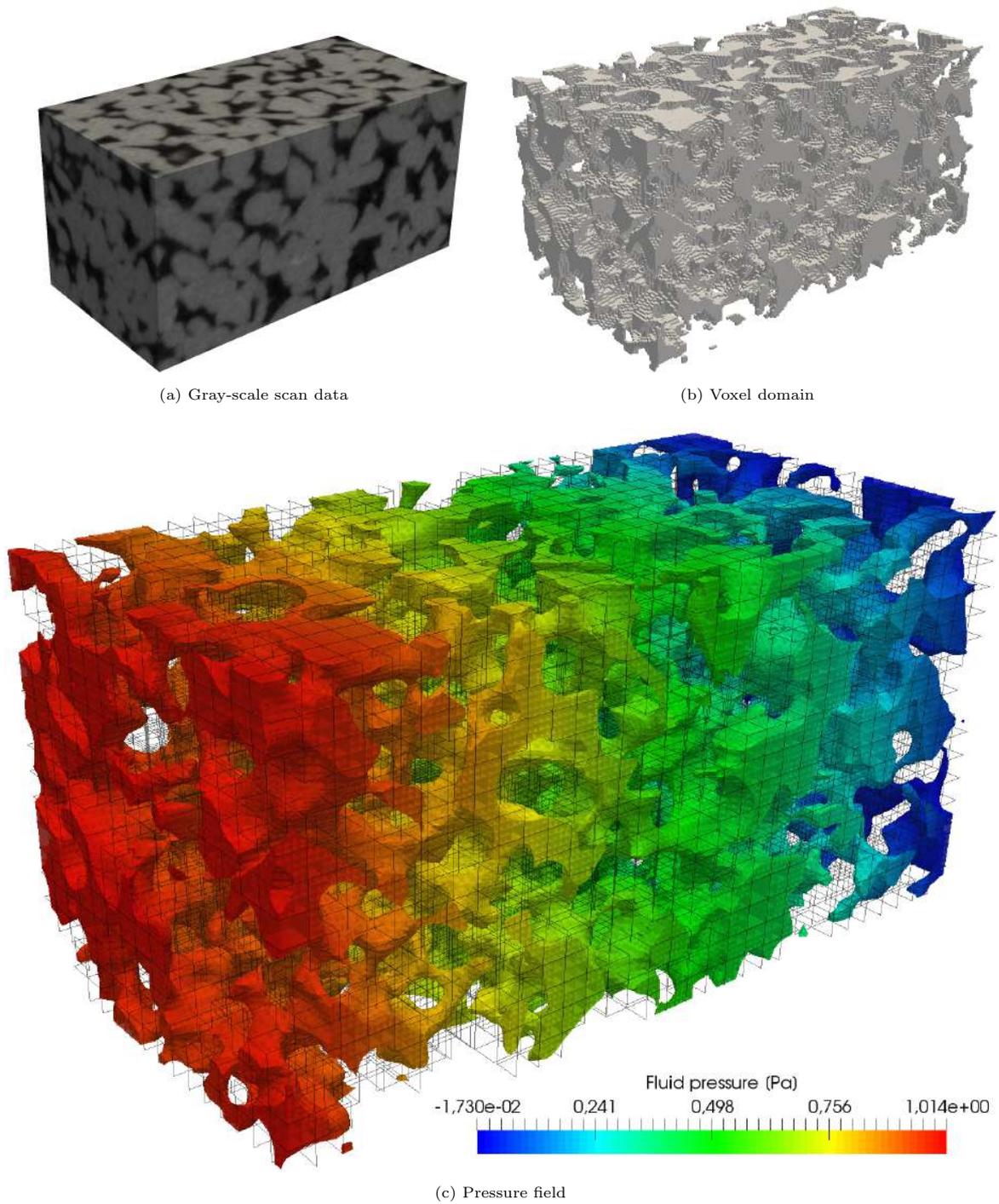

(a) Gray-scale scan data  (b) Voxel domain

(c) Pressure field

Figure 23: (a) $\mu$CT scan data consisting of $2 \times 10^6$ gray-scale voxels with a resolution of 25 $\mu$m. (b) Voxel mesh obtained by segmentation of the scan data. (c) Pressure field computed using skeleton-stabilized immersogeometric analysis on a $16 \times 16 \times 32$ background mesh with $\varrho_{\max} = 1$ and quadratic $C^1$ splines.

# 6. Conclusions

A skeleton-stabilized formulation was proposed for immersed isogeometric analyses of unsteady incompressible flow problems using identical B-spline bases for the velocity and pressure fields. This formulation extends the developments in Ref. [22], where mesh-conforming isogeometric analysis of incompressible flow problems using identical pressure and velocity bases was considered. The pivotal idea behind the considered stabilization technique – which can be regarded as the isogeometric extension of the (continuous) interior penalty method – is that the inf-sup condition is bypassed by supplementing a penalty term for basis function derivative jumps across element interfaces, thereby effectively penalizing oscillatory pressure behavior. The mesh-conforming formulation is amended with a ghost-penalty stabilization to resolve ill-conditioning issues and maintain the robustness in case of small volume fraction cut cells in the immersed finite cell setting. This ghost-penalty term bears close resemblance to the pressure-stabilization operator, but acts on the velocity field in the vicinity of the immersed boundaries only.

An important aspect of this work is that we fully leverage the smoothness properties of the full-regularity B-spline basis functions constructed over the background mesh, in the sense that the stabilization operators only act on the highest-order normal derivatives of the basis functions via their interface jumps. All lower order derivatives vanish as a result of the continuity properties of the B-spline basis. One advantage of this isogeometric approach in comparison to the case of a Lagrange-based analysis is that it only requires a penalty parameter for the highest-order derivative jump. Moreover, the impact of the stabilization term on the sparsity pattern of the system matrix is significantly reduced when full-regularity B-splines are considered instead of Lagrange finite elements. This is an important benefit of the considered isogeometric approach from the perspective of computational effort.

In a series of two-dimensional benchmarks problems we have observed optimal rates of convergence for the $L^2$ and $H^1$-norms of the velocity field and the $L^2$ norm of the pressure field. In comparison to the immersed simulation results based on inf-sup stable isogeometric finite element pairs considered in Ref. [21], we observed oscillation-free pressure fields near the cut boundaries. As a result, using the skeleton-based stabilization technique considered herein, quantities of interest pertaining to the immersed boundaries can be computed reliably. Provided that the immersed boundary is parametrized with sufficient smoothness using the employed bi-section based tessellation scheme, double-order convergence rates have been observed for (sufficiently smooth) quantities of interest pertaining to the immersed boundary.

To demonstrate the potential of the immersed isogeometric framework for geometrically and topologically complex three-dimensional domains, an image-based three-dimensional analysis of a flow through a porous medium was presented. For this type of geometry – which is defined by segmentation of the smoothened $\mu$CT-scan voxel data – virtually all elements are cut, which renders the standard cell-based techniques such as inf-sup stable space and Galerkin-Least square approaches ineffective. It was demonstrated that the proposed framework can yield results that are in good correspondence with finite volume based reference results on meshes that are significantly coarser than the voxel meshes used for the reference simulation. Thus the proposed approach alleviates the constraint between the resolution of the image data and that of the computational meshes.

The simulations considered herein were restricted to moderate Reynolds number flows on uniform background meshes. In this setting satisfactory results are obtained when the skeleton-penalty and ghost-penalty parameters are chosen in adequate ranges. The rates of convergence have been observed to be insensitive to the choice of the penalization parameters for a wide range of considered values. Consistent with observation in the conforming case, the magnitude of the $L^2$ pressure error is observed to be influenced by the choice of the parameters, but a significant range of parameters exists for which this error is insensitive to the precise values of the parameters. This makes ad hoc selection of the parameters practical. The development of more specific selection criteria – preferably in the form of rigorously derived explicit expressions – is an import aspect of the further development of the proposed simulation framework.

As part of the further development of the mathematical analysis of the proposed formulation, in our future work we aim at obtaining a more fundamental understanding of the influence of geometric irregularities on the approximation quality. In relation to this, we also aim at exploiting the locally refined spline discretizations that can be constructed over the regular background mesh, and to use these refinements in a mesh-adaptive analysis. Extension to convection-dominated problems – which is a non-trivial work in the sense that an additional convection-stabilization technique must be combined with the stabilization techniques already considered – is also an important topic of further study.

**Acknowledgements** We acknowledge the support from the European Commission EACEA Agency, Framework Partnership Agreement Ref. 2013-0043 Erasmus Mundus Action 1b, as a part of the *EM Joint Doctorate Simulation in Engineering and Entrepreneurship Development* (SEED). T.H. would like to thank the support of the European Research Council under the ERC-AG Agreement No. 341225. A.R. also acknowledges the support of Fondazione Cariplo - Regione Lombardia through the project "Verso nuovi strumenti di simulazione super veloci ed accurati basati sull'analisi isogeometrica", within the program RST - rafforzamento. C.Q. acknowledges the support of the Darcy Centre of Eindhoven University of Technology and Utrecht University. C.Q. and H.v.B. acknowledge the support of the Netherlands Organization for Scientific Research (NWO) Industrial Partnership Programme "Fundamental Fluid Dynamics Challenges in Inkjet Printing" (IPP FIP).

We thank Michael Afanasyev from Delft University of Technology for providing the 3D $\mu$CT data. The simulations in this work were performed using the open source software package Nutils ([51], www.nutils.org).